\documentclass[runningheads]{llncs}

\usepackage{tikz}
\usetikzlibrary{arrows.meta}
 \usetikzlibrary{decorations.pathreplacing}

\usepackage{amsmath,amssymb,mathtools}
\usepackage{booktabs}
\usepackage{enumitem}
\usepackage{microtype}
\usepackage[hidelinks]{hyperref}
\usepackage[nameinlink,noabbrev]{cleveref}
\usepackage{todonotes}
\crefname{assumption}{assumption}{assumptions}
\Crefname{assumption}{Assumption}{Assumptions}
\crefname{remark}{remark}{remarks}
\Crefname{remark}{Remark}{Remarks}
\usepackage{bm}

\usepackage[numbers,square]{natbib}

\newcommand{\cA}{\mathcal{A}}
\newcommand{\cB}{\mathcal{B}}
\newcommand{\cC}{\mathcal{C}}

\newcommand{\cI}{\mathcal{I}}

\newcommand{\cN}{\mathcal{N}}
\newcommand{\cP}{\mathcal{P}}

\newcommand{\cS}{\mathcal{S}}
\newcommand{\cT}{\mathcal{T}}

\newcommand{\cV}{\mathcal{V}}
\newcommand{\cX}{\mathcal{X}}
\newcommand{\cY}{\mathcal{Y}}
\newcommand{\cZ}{\mathcal{Z}} 

\newcommand{\lambdamin}{\lambda_{\min}}
\newcommand{\lambdamax}{\lambda_{\max}}

\newcommand{\R}{\mathbb{R}}
\newcommand{\Z}{\mathbb{Z}}

\newcommand{\Pp}{\mathbb{P}}

\newcommand{\one}{\mathbf{1}}
\newcommand{\zero}{\mathbf{0}}
\newcommand{\conv}{\operatorname{conv}}

\newtheorem{assumption}{Assumption}

\usepackage[T1]{fontenc}

\usepackage{newtxtext}       %
\usepackage[varvw]{newtxmath}       %

\usepackage{graphicx}
\begin{document}
\title{Probabilistically Exact Convexification of Mixed-Integer Quadratic Optimization with Indicators}
\titlerunning{Probabilistically Exact Convexification}
\author{Salar Fattahi\inst{1} \and
Andr{\'e}s G{\'o}mez\inst{2}}
\authorrunning{Fattahi and G{\'o}mez}
\institute{University of Michigan, Ann Arbor MI, USA
\email{fattahi@umich.edu}\\ 
\and
University of Southern California, Los Angeles CA, USA
\email{gomezand@usc.edu}}
\maketitle              %
\begin{abstract}
We study the convexification of mixed-integer quadratic optimization with indicator variables. The problem remains NP-hard even when the quadratic matrix has bandwidth two, posing a fundamental obstacle to tractable exact convexification. To overcome this challenge in realistic settings, we introduce \emph{probabilistically exact (inner) convexification}, a framework that relaxes the requirement of exactness for every objective while preserving exact optimality for {\it statistically relevant ones}. Our approach is motivated by inference applications, where the objective is determined by a data vector drawn from a distribution and a regularization vector chosen from a prescribed family. Our convexification preserves the exactness of the optimal value with high probability over the data distribution, simultaneously for every regularization vector in this family.
We apply this notion to quadratic matrices of any fixed bandwidth. Despite their worst-case NP-hardness, we construct probabilistically exact convexifications that admit exact extended second-order cone formulations of size polynomial in the dimension $n$ and the inverse failure probability $1/\delta$. The required assumptions for their high-probability exactness are surprisingly mild: they accommodate both light-tailed and heavy-tailed data, including distributions without finite mean, as well as arbitrary nonnegative regularization coefficients.

\end{abstract}
\section{Introduction}
Mixed-integer quadratic optimization with indicator variables provides a natural framework for statistical estimation with logical decisions. Given a positive definite matrix $\bm Q\in\R^{n\times n}$, we consider optimization problems of the form
\begin{equation}\label{eq:MIQP}
    \min_{(t,\bm x,\bm z)\in\cP}
    \left\{t+\bm\alpha^\top\bm x+\bm\beta^\top\bm z\right\},
\end{equation}
\begin{equation}\label{eq::P}
    \text{where}\quad \cP:=\left\{(t,\bm x,\bm z)\in\R_+\times\R^n\times\{0,1\}^n:
    t\geq\bm x^\top\bm Q\bm x,\quad
    \bm x\circ(\one-\bm z)=\zero\right\},
\end{equation}
and $\circ$ denotes the entrywise product. The indicator constraints  $\bm x\circ(\one-\bm z)=\zero$ are equivalent to the logical condition $x_i\neq 0\implies z_i\neq 0$ for all $i\in [n]$; such constraints are used to impose logical conditions on the support of the continuous variables $\bm{x}$. 

An application particularly relevant to this paper is sparse inference of spatio-temporal processes~\cite{besag1991bayesian}. Given noisy observations $\bm\alpha\in\R^n$ and a sparse positive definite  matrix $\bm Q\in\R^{n\times n}$, a representative formulation is~\cite{atamturk2021sparse,gomez2021outlier,gomez2025convex}
\begin{align}\label{eq:BYM}
    \min_{\substack{\bm x\in \mathbb R^n,\bm z\in\{0,1\}^n}}\quad
    &\frac{1}{2}\|\bm\alpha-\bm x\|_2^2+\bm x^\top\bm Q\bm x
    +\beta\one^\top\bm z \nonumber\\
    \text{s.t.}\quad
    &\bm x\circ(\one-\bm z)=\zero.
\end{align}
Variables $\bm x$ represent estimates of the underlying process, while the indicators $z_i$ determine which signal components $x_i$ may be nonzero. The parameter $\beta\geq 0$ penalizes the number of active components, with $\one^\top\bm z$ counting the nonzero signal components at optimality. Finally, $\bm x^\top\bm Q\bm x$ regularizes the signal according to its dependence structure, where $Q_{ij}\neq 0$ indicates that components $i$ and $j$ are spatially or temporally related. In particular, $\bm Q$ is typically banded in time-series models~\cite{liu2023graph,gomez2024real}.

Convexification-based approaches to solving~\eqref{eq:MIQP} exploit convex formulations of the underlying mixed-integer set~\cite{wei2022ideal,wei2024convex}. An exact characterization of the convex hull $\operatorname{conv}(\cP)$ yields an equivalent convex reformulation of~\eqref{eq:MIQP}. For a fixed quadratic matrix $\bm Q$, this formulation can be reused across different (data) vectors $\bm{\alpha}$ and (hyperparameter) vectors $\bm{\beta}$. When tractable, it therefore facilitates repeated optimization and cross-validation in statistical applications.
However, obtaining a tractable exact convexification of $\cP$ is challenging, and available results rely on additional structural assumptions. For example, tridiagonal quadratics admit exact formulations of quadratic size~\cite{liu2023graph,gomez2024real}, while decision diagrams provide a framework for constructing exact and approximate formulations for broader classes of structured quadratics~\cite{choi2026convexification}. Despite these advances, tractable exact convexification remains elusive in general: problem~\eqref{eq:MIQP} remains NP-hard even when $\bm Q$ is banded with bandwidth two~\citep[Theorem~1 and Corollary~2]{m18_bandwidth_two_2026}.\footnote{The cited result appears in an AI-generated note released in the repository accompanying \cite{gusev2026swarms}. {While the authors from the repository do not claim to have verified correctness of this note, } we have independently checked the proof of the cited result and verified its correctness.} Consequently, banded $\bm Q$ alone does not guarantee an exact convex reformulation that can be constructed and solved efficiently.

To alleviate the burdens associated with exact convexifications, which guarantee exactness for {\it all} possible data vectors and regularization parameters $(\bm{\alpha},\bm{\beta})$, we ask whether it is possible to construct simpler formulations that are exact only for the \emph{relevant statistical inputs}. Indeed, in statistical applications, the regularization parameters $\bm{\beta}$ are typically nonnegative. Moreover, the data vector $\bm{\alpha}$ is often modeled as a random draw from a well-behaved distribution. Motivated by these observations, we show that, when $\bm{Q}$ is banded, it is possible to construct compact (inner) convex approximations of the closure of $\operatorname{conv}(\mathcal{P})$ that recover exact optimal solutions of~\eqref{eq:MIQP} {\it with high probability} over the distribution of $\bm\alpha$, uniformly for all nonnegative choices of $\bm \beta$.

Due to space constraints, most proofs are deferred to the appendix.

\section{Probabilistically Exact Approximation and Convexification}
For any closed set
\(\cC\subseteq\mathbb{R}_{+}\times\mathbb{R}^n\times\mathbb{R}^n\),
define its \textit{restricted support function}
\(h_{\cC}:\mathbb{R}^n\times\mathbb{R}^n\to\overline{\mathbb{R}}\) by
\begin{align}
    h_{\cC}(\bm{\alpha},\bm{\beta})
    :=
    \inf_{(t,\bm{x},\bm{z})\in\cC}
    \left\{
        t+\bm{\alpha}^\top \bm{x}+\bm{\beta}^\top \bm{z}
    \right\},
\end{align}
where
\(\overline{\mathbb{R}}:=\mathbb{R}\cup\{-\infty,+\infty\}\), and we adopt
the convention that \(h_{\cC}(\bm{\alpha},\bm{\beta})=+\infty\) when
\(\cC=\varnothing\).
The term ``restricted'' refers to the fact that \(h_{\cC}\) is the
restriction of the full support function
to directions for which the coefficient of \(t\) is equal to one. 
\begin{definition}[Upward-closedness]
A set $\cC\subseteq\mathbb{R}_+\times \mathbb{R}^n\times\mathbb{R}^n$
is said to be \emph{upward-closed (in its first element)} if $\cC$ is closed and
\[
    (t, \bm{x},\bm{z})\in \cC
    \quad\Longrightarrow\quad
    (t+s, \bm{x},\bm{z})\in \cC
    \qquad \text{for every }s\geq 0.
\]
\end{definition}
\begin{lemma}[Convex hull via restricted support functions]\label{lem:restricted-support-characterization}
    Suppose that
    \(\cX\subseteq\cP \) is nonempty and
    upward-closed. Then, both $\conv(\cX)$ and $\conv(\cP)$ are closed and
    \begin{align}
        h_{\cX}(\bm{\alpha},\bm{\beta})
        =
        h_{\cP }(\bm{\alpha},\bm{\beta}),
        \quad
        \forall(\bm{\alpha},\bm{\beta})\in\mathbb{R}^n\times\mathbb{R}^n, \quad \iff \quad \conv(\cX)=\conv(\cP ).
    \end{align}
\end{lemma}
The preceding lemma characterizes the exact convexification of $\cP$ by equality of the restricted support functions in {all} directions $(\bm{\alpha}, \bm{\beta})$.
Motivated by this characterization, we introduce a relaxed notion of convexification for \(\cP\): rather than requiring the restricted support function of \(\cX\subseteq\cP\) to coincide with that of \(\cP\) for every pair \((\bm{\alpha},\bm{\beta})\), we require exactness only for \(\bm{\beta}\) belonging to a prescribed set \(\Gamma\subseteq\mathbb{R}^n\) and for a set of coefficients \(\bm{\alpha}\) having sufficiently high probability under a probability measure \(\mu\).
\begin{definition}[probabilistically exact approximation and convexification]\label{def:inner-approx-convex}
    Let $\Gamma\subseteq\mathbb{R}^n$, let $\mu$ be a probability measure on $\mathbb{R}^n$, and let
    $\delta\in[0,1]$. A set
    $\cX\subseteq\cP $ is called a
    \emph{$(\mu,\Gamma,\delta)$-probabilistically exact (inner) approximation}
    of $\cP $ if it is nonempty and upward-closed, and satisfies
    \begin{align}\label{eq::prob-condition}
        \mathbb{P}_{\bm\alpha\sim\mu}
        \left(
            h_{\cX}(\alpha,\beta)
            =
            h_{\cP }(\alpha,\beta)
            \quad\text{for all }\beta\in\Gamma
        \right)
        \geq 1-\delta.
    \end{align}
    The corresponding convex set $\conv(\cX)$ is called a \emph{$(\mu,\Gamma,\delta)$-probabilistically exact (inner) convexification}
    of $\cP $.
\end{definition}
The definition is tailored to realistic data-driven settings in which neither \(\bm{\alpha}\) nor \(\bm{\beta}\) is treated as arbitrary. In particular, inspired by the inference problem~\eqref{eq:BYM}, the (data) vector \(\bm{\alpha}\) defining the linear objective in the \(x\)-variables is modeled as a random vector with distribution \(\mu\), whereas the (regularization) vector \(\bm{\beta}\) is selected from a prescribed set \(\Gamma\) (e.g., nonnegative vectors). The condition \eqref{eq::prob-condition} requires that, with probability at least \(1-\delta\) over \(\bm{\alpha}\sim\mu\), restricting the feasible region from \(\cP\) to \(\cX\) preserves the optimal value simultaneously for every \(\bm{\beta}\in\Gamma\). Thus, \(\cX\) need only retain the portions of \(\cP\) relevant to objective coefficients that occur with sufficiently high probability under \(\mu\); exactness may fail on a set of coefficients \(\bm{\alpha}\) having probability at most \(\delta\). 
The key advantage of probabilistic approximation and convexification is that they can yield more tractable formulations by allowing a controlled probability of inexactness. We demonstrate this advantage for the important class of problems with banded $\bm Q$: although
linear optimization over $\cP$ remains NP-hard for general banded
matrices, we construct a compact probabilistically exact
convexification under mild assumptions.

\section{Application to Banded Structures}
In this section, we derive polynomial-size probabilistically exact approximation and convexification of $\cP$ when $\bm Q$ is {banded}. For an integer $b\ge1$, we call $\bm Q$ \textit{$b$-banded} if $\bm Q_{ij}=0$ whenever $|i-j|>b$. 
\begin{assumption}[Bounded spectrum and nonnegative regularization]\label{ass:deterministic}
There are constants $0<\lambda_{\min}\le\lambda_{\max}<\infty$ and $\beta_0>0$, independent of $n$, such that
\begin{align}
\lambdamin \bm I&\preceq \bm Q\preceq\lambdamax \bm I, \label{eq:spectrum}\\
\Gamma&\subseteq\{\zero\}\cup[\beta_0,\infty)^n. \label{eq:Gamma}
\end{align}
\end{assumption}
Assumption~\eqref{eq:Gamma} requires positive regularization parameters to be bounded away from zero, consistent with positive penalties used to promote sparsity in $\bm x$. We also allow $\bm\beta=\zero$ to accommodate the unregularized problem, which is commonly considered in cross-validation. 
We summarize our main result in the next informal theorem.

\begin{theorem}[Probabilistic exactness for banded structures; informal]\label{thm::banded}
Suppose that $\bm Q$ is $b$-banded and satisfies~\eqref{eq:spectrum}, $\Gamma$ satisfies~\eqref{eq:Gamma}, and $\mu$ is a product measure whose coordinates satisfy: $(i)$ uniformly positive probability near zero; and $(ii)$ uniformly bounded absolute moments of some positive order.
Then, for any $\delta\in(0,1]$, there exists a $(\mu,\Gamma,\delta)$-probabilistically exact approximation $\cX$ of $\cP$ such that $\conv(\cX)$ is closed and admits an exact extended second-order cone programming (SOCP) formulation with
\begin{equation}\label{eq:ideal-DD-size}
    \begin{cases}
        O\!\left(
            n\left[
                1+\log\left(\dfrac{n}{\delta}\right)
            \right]^2
        \right),
        & b=1,\\[3mm]
        O\!\left(\vartheta_b^b
            n\left[
                1+\log\left(\dfrac{n}{\delta}\right)
            \right]\left(\dfrac{n}{\delta}\right)^{\gamma_b}
        \right),
        & b\ge2,
    \end{cases}
\end{equation}
variables and constraints. Here, for $b\ge2$, $\vartheta_b\in(1,2)$ is the unique positive root of
$x^b=x^{b-1}+\cdots+x+1$, and $\gamma_b=O(\log\vartheta_b)$.
In particular, for any fixed bandwidth $b$, the formulation has size polynomial in $n$ and $1/\delta$.
\end{theorem}
The assumptions on $\Gamma$ and $\mu$ are satisfied in a broad range of practical settings. In particular, the restrictions on $\Gamma$ only exclude negative regularization parameters, consistent with sparsity-inducing regularization in practice. Similarly, the conditions on $\mu$ accommodate a broad class of distributions, including Gaussian, Laplace, and bounded distributions with sufficient probability mass near zero, as well as heavy-tailed distributions such as Cauchy and Student-$t$. Notably, the latter class includes distributions with infinite mean.

\paragraph{Comparison with existing results.}
For tridiagonal matrices ($b=1$), \citet{liu2023graph} establish an exact semidefinite programming formulation of $\conv(\cP)$ with $O(n^2)$ variables and linear constraints, together with a single $(n+1)\times(n+1)$ positive semidefinite constraint. Subsequently, \citet{gomez2024real} derive an exact SOCP formulation with $O(n^2)$ variables and constraints. Under our distributional assumptions, we improve upon these size bounds by constructing a probabilistically exact convexification with only $O(n\log^2 n)$ variables and constraints for fixed $\delta$ and model parameters.
For general banded matrices with $b\ge2$, the situation is fundamentally different: problem~\eqref{eq:MIQP} remains NP-hard even for $b=2$~\citep[Theorem~1 and Corollary~2]{m18_bandwidth_two_2026}, implying that bounded bandwidth alone cannot guarantee an efficiently constructible and solvable exact convex formulation unless $\mathrm{P}=\mathrm{NP}$. To address this challenge, \citet{gomez2024real} and \citet{bienstock2024solving} develop a deterministic approximation scheme that obtains the optimal value within an additive error $\epsilon$, with complexity polynomial in $n$ and $\|\bm\alpha\|_\infty^2/\epsilon$ for fixed bandwidth and spectral parameters. Our approach offers a different tradeoff between tractability and accuracy: rather than allowing a bounded approximation error for every instance, we construct polynomial-size convexifications that preserve the optimal value \emph{exactly} with probability at least $1-\delta$ over $\bm\alpha\sim\mu$. Moreover, this exactness guarantee holds uniformly for all $\bm\beta\in\Gamma$, allowing the same formulation to be reused across different regularization choices without sacrificing optimality on the high-probability event.

\subsection{Key Idea: Dynamic Programming}\label{subsec::DP}
Our construction relies on dynamic programming (DP). Throughout the paper, we index the variables according to the row and column ordering of $\bm Q$, and write $[\ell:k]:=\{\ell,\ell+1,\ldots,k\}$ for integers $\ell\le k$.
A support vector $\bm z\in\{0,1\}^n$ is said to contain a \emph{resetting sequence of length $b$} if there exists $s\in[1:n-b+1]$ such that $z_s=\cdots=z_{s+b-1}=0$. Since $\bm Q$ is $b$-banded, these $b$ consecutive zeros eliminate all quadratic interactions between variables indexed before and after the resetting sequence. This separation allows the problem to be decomposed into smaller subproblems, forming the basis of our DP construction.

\begin{figure}[t]
\centering
\begin{tikzpicture}[
  x=1.03cm, y=1cm, font=\small,
  vertex/.style={
    circle, draw=black, line width=.45pt, fill=white,
    minimum size=5.5mm, inner sep=0pt, font=\footnotesize
  },
  selected/.style={vertex, fill=black!25},
  removed/.style={vertex, draw=black!45, text=black!50},
  edge/.style={draw=black!75, line width=.45pt},
  removed edge/.style={draw=black!50, dashed, line width=.45pt},
  group/.style={
    decorate, decoration={brace, mirror, amplitude=4pt},
    line width=.45pt
  },
  note/.style={font=\small, align=center, inner sep=1pt}
]

\foreach \i in {1,...,12} {
  \coordinate (b\i) at ({\i-1},0);
}

\path[
  fill=black!5, draw=black!50, dashed,
  rounded corners=3pt, line width=.45pt
] (3.62,-0.72) rectangle (5.38,0.72);

\foreach \i/\j in {1/2,2/3,3/4} {
  \draw[edge] (b\i) -- (b\j);
}
\foreach \i/\j in {1/3,2/4} {
  \draw[edge] (b\i) to[out=42,in=138] (b\j);
}

\foreach \i/\j in {7/8,8/9,9/10,10/11,11/12} {
  \draw[edge] (b\i) -- (b\j);
}
\foreach \i/\j in {7/9,8/10,9/11,10/12} {
  \draw[edge] (b\i) to[out=42,in=138] (b\j);
}

\foreach \i/\j in {4/5,5/6,6/7} {
  \draw[removed edge] (b\i) -- (b\j);
}
\foreach \i/\j in {3/5,4/6,5/7,6/8} {
  \draw[removed edge] (b\i) to[out=42,in=138] (b\j);
}

\foreach \i in {1,3,4,7,9,10,12} {
  \node[selected] at (b\i) {$\i$};
}

\foreach \i in {2,5,6,8,11} {
  \node[removed] at (b\i) {$\i$};
  \begin{scope}[shift={(b\i)}]
    \draw[black!55, line width=.55pt]
      (-.22,-.22) -- (.22,.22);
  \end{scope}
}

\draw[group] (-.28,-0.77) -- (3.28,-0.77);
\node[note, anchor=north, text width=3.9cm] at (1.5,-1)
  {Left subproblem $f_4^\star$};

\draw[group] (3.72,-0.77) -- (5.28,-0.77);
\node[note, anchor=north, text width=2.45cm] at (4.5,-1)
  {Resetting sequence};

\draw[group] (5.72,-0.7) -- (11.28,-0.7);
\node[note, anchor=north, text width=5.0cm] at (8.5,-1)
  {Reset-free subproblem $q^\star_{[7:12]}$};

\end{tikzpicture}
\caption{Illustration of the DP recurrence for $b=2$ and $n=12$. The last resetting sequence, $z_5^\star=z_6^\star=0$, separates the indices $[1:4]$ from the reset-free indices $[7:12]$.}
\label{fig:dp-recurrence-b2}
\end{figure}
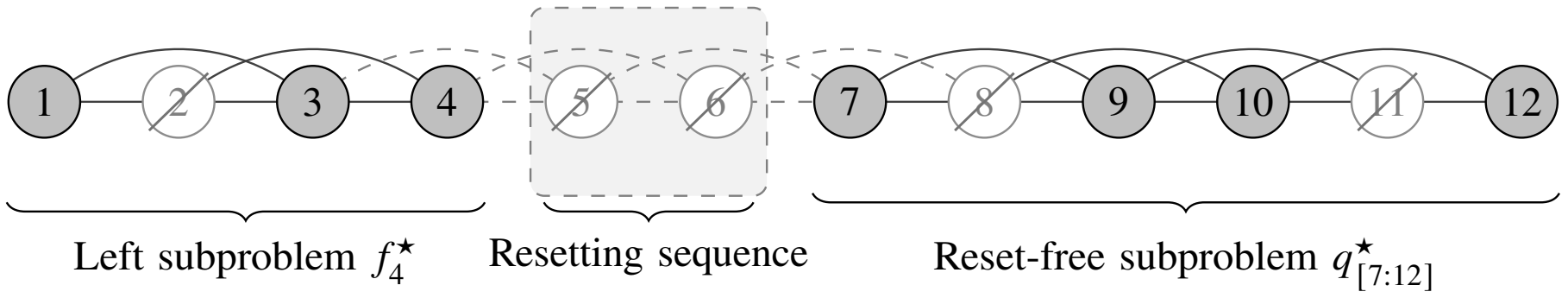

To formalize this decomposition, we introduce two families of subproblems. For $s\in\{0,\ldots,n\}$, let $f_s^\star$ denote the optimal value of~\eqref{eq:MIQP} with $z_{s+1}=\cdots=z_n=0$. For $1\le\ell\le k\le n$, let $q_{[\ell:k]}^\star$ denote the optimal value over support vectors restricted to $[\ell:k]$ and containing no resetting sequence entirely within this interval. Formally,
\begin{align}\label{eq:reset-free}
    q^\star_{[\ell:k]}
    = \min_{\substack{(t,\bm x,\bm z)\in\cP\\
    \bm z\in\cN_{[\ell:k],b}}}
    \left\{t+\bm\alpha^\top\bm x+\bm\beta^\top\bm z\right\},
\end{align}
where
\begin{align}\label{eq:reset-free-support}
    \cN_{[\ell:k],b}
    :=
    \left\{
    \bm z\in\{0,1\}^n:
    \begin{array}{ll}
        z_i=0, & \forall i\notin[\ell:k],\\
        \displaystyle\sum_{j=s}^{s+b-1}z_j\ge1,
        & \forall s\in[\ell:k-b+1]
    \end{array}
    \right\}.
\end{align}
Let $\bm z^\star\in\{0,1\}^n$ be an optimal support vector for $f^\star_i$, and let $s$ denote the starting index of its last resetting sequence of length $b$, with $s\leq 0$ if no such sequence exists. Conditioned on the location of $s$, we obtain the DP recurrence
\begin{align}
    f_i^\star
    =
    \min_{-b+1\le s\le i-b+1}
    \left\{
        f_{s-1}^\star+q_{[s+b:i]}^\star
    \right\},
    \qquad i\in[n],
    \label{eq:DP-recurrence}
\end{align}
where $f_j^\star=0$ for $j\le0$ and $q_{[i+1:i]}^\star=0$.
In particular, $f_n^\star$ is the optimal value of~\eqref{eq:MIQP}; Figure~\ref{fig:dp-recurrence-b2} exemplifies this decomposition. 
Evidently, the efficiency of this DP hinges on the cost of computing the values $q_{[\ell:k]}^\star$. Once these are computed, the partial costs $f_1^\star,\ldots,f_n^\star$ can be obtained sequentially using~\eqref{eq:DP-recurrence}.
We distinguish two cases:

\begin{itemize}
\item \textbf{Tridiagonal case ($\bm{b=1}$).}
For $b=1$, $\cN_{[\ell:k],1}$ is a singleton, with $z_i=1$ for $i\in[\ell:k]$ and $z_i=0$ otherwise. Thus, by \eqref{eq:reset-free}, each $q_{[\ell:k]}^\star$ can be computed by solving a single QP. In fact, \citet{liu2023graph} show that all $\left\{q_{[\ell:k]}^\star\right\}_{\ell\le k}$ can be computed in $O(n^2)$ time. Consequently,~\eqref{eq:DP-recurrence} yields an $O(n^2)$-time algorithm for solving~\eqref{eq:MIQP}.

\item \textbf{General banded case ($\bm{b\ge2}$).}
For $b\ge2$, we show that
$|\cN_{[\ell:k],b}|=\Theta\left(\vartheta_b^{k-\ell+1}\right)$,
where $\vartheta_b\in(1,2)$ is the unique positive root of
\[
    x^b=x^{b-1}+\cdots+x+1.
\]
Hence, directly evaluating each $q_{[\ell:k]}^\star$ by enumerating admissible supports in $\cN_{[\ell:k],b}$ can require exponentially many QP solves, making this approach computationally prohibitive. This difficulty is expected in light of the NP-hardness of~\eqref{eq:MIQP} for $b\geq 2$.
\end{itemize}
For $b\ge2$, the computational bottleneck in~\eqref{eq:DP-recurrence} is evaluating $q^\star_{[s+b:i]}$ over long intervals. Our key idea is therefore to avoid these costly subproblems by restricting the feasible region to support vectors whose reset-free intervals have length at most a prescribed threshold. 

\begin{definition}[Reset support sets]
\label{def:reset-support-set}
For integers $b\in [n]$ and \(M\geq b\), define the \emph{\((M,b)\)-reset support set} by
\begin{equation}\label{eq:reset-language}
    \cZ_{M,b}
    :=
    \left\{
        \bm z\in\{0,1\}^n :
        \begin{array}{l}
            \text{every block of $M$ consecutive entries of $\bm z$}\\
            \text{contains a resetting sequence of length $b$}
        \end{array}
    \right\}.
\end{equation}
When $M>n$, the condition is vacuous.
Additionally, define the
corresponding \emph{augmented \((M,b)\)-reset support set} by
\begin{equation}\label{eq:augmented-language}
    \widehat{\cZ}_{M,b}
    :=
    \begin{cases}
        \cZ_{M,b}\cup\{\one\}, & \zero\in\Gamma,\\
        \cZ_{M,b}, & \zero\notin\Gamma.
    \end{cases}
\end{equation}
\end{definition}
We define the associated restricted feasible set by
\begin{equation}\label{eq:restricted-P}
    \cP_{M,b}
    :=
    \left\{
        (t,\bm x,\bm z)\in\cP :
        \bm z\in\widehat{\cZ}_{M,b}
    \right\}.
\end{equation}
Apart from the possible inclusion of \(\bm z=\one\), the set
\(\cP_{M,b}\) restricts \(\cP\) to support vectors for which every
block of \(M\) consecutive coordinates contains a resetting sequence
of length \(b\). To formulate the corresponding DP, let $\widehat f_{i}$ and
$\widehat q_{[\ell:k]}$ be defined as $f_i^\star$ and
$q_{[\ell:k]}^\star$, respectively, with the additional constraint
$\bm z\in\cZ_{M,b}$. 
The resulting DP recurrence is
\begin{align}
    \widehat f_{i}
    =
    \min_{\substack{1-b\le s\le i-b+1\\
                   |[s+b:i]|\le M-1}}
    \left\{
        \widehat f_{s-1}
        +\widehat q_{[s+b:i]}
    \right\},
    \qquad i\in[n],
    \label{eq:DP-recurrence-restricted}
\end{align}
where $\widehat f_{j}=0$ for $j\le0$,
$\widehat q_{[i+1:i]}=0$, and infeasible subproblems have value
$+\infty$.
Unlike~\eqref{eq:DP-recurrence}, this restricted recurrence considers only reset-free intervals of length at most $M-1$. Each interval cost $\widehat q_{[\ell:k]}$ appearing in~\eqref{eq:DP-recurrence-restricted} can therefore be computed by solving at most $O(\vartheta_b^M)$ QPs. Since~\eqref{eq:DP-recurrence-restricted} involves $O(nM)$ such costs, computing all of them requires $O(nM\vartheta_b^M)$ QP solves. Once these costs are available, the recurrence computes $\widehat f_n$ in $O(nM)$ time. If $\zero\in\Gamma$, we also evaluate the full support $\bm z=\one$ by solving one additional QP and take the smaller of its optimal value and $\widehat f_n$. This yields the optimal value of~\eqref{eq:MIQP} over the restricted feasible set $\cP_{M,b}$. This approach raises two natural questions:

\paragraph{1. Does restricting $\cP$ to $\cP_{M,b}$ preserve optimality?}
In Section~\ref{subsec::barriers}, we show that, if $\Gamma$ satisfies~\eqref{eq:Gamma} and $\mu$ satisfies a \emph{barrier-free property} (Definition~\ref{def:void-property}), there exists a choice of
\[
    M=O\!\left(\log\left(\frac{n}{\delta}\right)+b\right)
\]
for which $\cP_{M,b}$ and $\conv(\cP_{M,b})$ are, respectively, a $(\mu,\Gamma,\delta)$-probabilistically exact approximation and convexification of $\cP$. Thus, with probability at least $1-\delta$ over $\bm\alpha\sim\mu$, replacing $\cP$ with either set preserves the optimal value of~\eqref{eq:MIQP} simultaneously for every $\bm\beta\in\Gamma$. Moreover, in Section~\ref{subsec::barrier-free}, we show that the barrier-free property required for $\mu$ holds under mild and verifiable conditions, including those delineated in Theorem~\ref{thm::banded}.

\paragraph{2. Does restricting $\cP$ to $\cP_{M,b}$ guarantee tractability?}
As discussed above, for $b\geq 2$, solving~\eqref{eq:MIQP} over $\cP_{M,b}$ via the restricted DP recurrence~\eqref{eq:DP-recurrence-restricted} requires solving $O(nM\vartheta_b^M)$ QPs. 
More generally, in Section~\ref{subsec::SOCP}, we show that $\conv(\cP_{M,b})$ admits an exact extended SOCP formulation with $O(nM\Phi_b(M))$ variables and constraints, where
\[
    \Phi_b(M)
    :=
    \begin{cases}
        M, & b=1,\\
        \vartheta_b^M, & b\ge2.
    \end{cases}
\]
With $M=O(\log(n/\delta)+b)$, both the running time and the formulation size are polynomial in $n$ and $1/\delta$ for fixed bandwidth $b$ and model parameters.

\subsection{Barriers and Probabilistic Reset Support Sets}\label{subsec::barriers}
Under Assumption~\ref{ass:deterministic}, let
\begin{equation}\label{eq:decay-parameters}
\kappa:=\frac{\lambdamax}{\lambdamin},\qquad
\rho:=\left(\frac{\sqrt\kappa-1}{\sqrt\kappa+1}\right)^{1/b}\in[0,1),\qquad
\tau_b:=\lambdamin \rho^{b-1}\sqrt{\frac{\beta_0}{2\lambdamax}}.
\end{equation}
We focus on the nontrivial case $\rho\in(0,1)$. If $\rho=0$, then $\kappa=1$ and $\bm Q$ is a scalar multiple of the identity, in which case $\operatorname{conv}(\cP)$ admits the standard linear-size perspective formulation \cite{frangioni2006perspective, gunluk2010perspective}.
We begin by introducing a local weight whose value triggers a reset.
\begin{definition}[Local weight]\label{def:local-weight}
    The \emph{local weight} at $i\in[n]$ is defined as $W_i(\bm\alpha) := \sum_{j=1}^n \rho^{|i-j|}|\alpha_j|$.
    We call $i\in[n-b+1]$ a \emph{barrier} when $W_i(\bm\alpha)\le\tau_b$.
\end{definition}
\begin{lemma}[Barriers force resets]
\label{lem:barrier-forces-zeros}
Fix $\bm\alpha\in\R^n$ and $i\in[n-b+1]$. If $W_i(\bm\alpha)\le\tau_b$, then, for every $\bm\beta\in[\beta_0,\infty)^n$, every optimal support vector $\bm z^\star$ of \eqref{eq:MIQP} has a resetting sequence of length $b$ beginning at $i$; that is, $z_i^\star=z_{i+1}^\star=\cdots=z_{i+b-1}^\star=0$.
\end{lemma}
By the above lemma, to ensure that every optimal solution belongs to \(\cP_{M,b}\), it suffices that, within every block of \(M\) consecutive coordinates, there exists at least one barrier \(i\) such that \(W_i(\bm\alpha)\le \tau_b\). We show next that this property holds with high probability whenever the probability measure \(\mu\) governing \(\bm\alpha\) satisfies a barrier-free property.

\begin{definition}[Barrier-free property]
\label{def:void-property}
Let \(C_{\rm v}\ge 1\) and \(c_{\rm v}>0\) be constants independent of \(n\).
A probability measure \(\mu\) on \(\R^n\) is said to satisfy the
\emph{\((C_{\rm v},c_{\rm v})\)-barrier-free property} if, for every interval
\(\cI\subseteq[n-b+1]\) of consecutive indices,
\begin{equation}\label{eq:void-property}
    \Pp_{\bm\alpha\sim\mu}
    \left(
        W_i(\bm\alpha)>\tau_b
        \ \text{for all } i\in\cI
    \right)
    \le
    C_{\rm v}\exp\!\left(-c_{\rm v}|\cI|\right).
\end{equation}
\end{definition}
\begin{theorem}[Probabilistic exactness under the barrier-free property]\label{thm:abstract-main}
Suppose that \(\mu\) satisfies the \((C_{\rm v},c_{\rm v})\)-barrier-free property. For any \(\delta\in(0,1]\), define
\begin{equation}\label{eq:L-M}
    M
    :=
    \left\lceil
        \frac{1}{c_{\rm v}}
        \log\left(\frac{C_{\rm v}n}{\delta}\right)
    \right\rceil
    +b-1.
\end{equation}
Then $\cP_{M,b}$ is a $(\mu,\Gamma,\delta)$-probabilistically exact approximation of $\cP$, and $\conv(\cP_{M,b})$ is a $(\mu,\Gamma,\delta)$-probabilistically exact convexification of $\cP$.
\end{theorem}
The preceding theorem establishes the probabilistic exactness of $\cP_{M,b}$ and its convex hull under the barrier-free property. To demonstrate the practical relevance of this result, we next show that the barrier-free property holds for broad classes of distributions.

\subsection{Measures Satisfying the Barrier-free Property}\label{subsec::barrier-free}
Our next result establishes the barrier-free property for product measures under two very mild conditions.

\begin{proposition}[Product measures satisfying the barrier-free property]
\label{prop:explicit-product-barrier-free}
Suppose that $\mu$ is a product measure on $\mathbb{R}^n$, that is, the coordinates of
$\bm\alpha\sim \mu$ are independent. Fix any $g\in(1,\rho^{-1})$ and define
\begin{equation}
    a_g
    :=
    \tau_b\frac{1-\rho g}{1+\rho g}
    >0,
\end{equation}
where $\tau_b$ is the barrier threshold defined in~\eqref{eq:decay-parameters}.
Assume that there exist constants $p_0\in(0,1]$,
$q>0$, and $M_q<\infty$, independent of $n$, such that
\begin{equation}
\label{eq:uniform-product-conditions}
    \min_{i\in[n]}
    \mathbb{P}(|\alpha_i|\leq a_g)
    \geq p_0,
    \qquad
    \max_{i\in[n]}
    \mathbb{E}|\alpha_i|^q
    \leq M_q.
\end{equation}
Then $\mu$ satisfies the $(C_v,c_v)$-barrier-free property for some constants
$C_v\geq 1$ and $c_v>0$ that depend only on
$\rho$, $g$, $\tau_b$, $p_0$, $q$, and $M_q$, and are independent of $n$.
\end{proposition}
In Appendix~\ref{app:explicit-dist}, we verify that these assumptions hold for independent coordinates with Gaussian, Laplace, Cauchy, or Student-$t$ distributions, as well as uniformly bounded coordinates satisfying the first condition in~\eqref{eq:uniform-product-conditions}. For example, for i.i.d.\ standard normal coordinates, one may choose $p_0=2\Phi(a_g)-1$, $q=1$, and $M_1=\sqrt{2/\pi}$, where $\Phi$ denotes the standard normal cumulative distribution function.

We devote the remainder of this subsection to outlining the proof of Proposition~\ref{prop:explicit-product-barrier-free}. A key idea underlying our analysis is to interpret the random coefficients $\alpha_1,\ldots,\alpha_n$ as inducing a random covering of the indices. For notational convenience, extend $\alpha_1,\ldots,\alpha_n$ to a sequence including large or negative indexes $(\alpha_j)_{j\in\mathbb{Z}}$, by setting $\alpha_j=0$ for every $j\notin[n]$. The added coordinates are deterministic and hence independent of the original coordinates, and they trivially satisfy the bounds in~\eqref{eq:uniform-product-conditions}. This extension allows us to formulate our argument without separately treating boundary effects.

\begin{definition}[Index coverage]
For each index $j\in\mathbb{Z}$, define its \emph{coverage radius} $R_j$ and \emph{coverage neighborhood} $\cB_j$ as 
\begin{equation}\label{eq:radius}
\cB_j
    :=
    \{i\in\mathbb{Z}: |i-j|<R_j\}, \quad\text{where}\quad 
    R_j
    :=
    \min\{k\in\mathbb{Z}_+ : |\alpha_j|\leq a_g g^k\}. 
\end{equation}
We say that an index $i\in\mathbb{Z}$ is \emph{covered} by
$j\in\mathbb{Z}$ if $i\in \cB_j$. Accordingly, an index $i\in\mathbb{Z}$ is
called \emph{uncovered} if it is not covered by any $j\in\mathbb{Z}$.
\end{definition}

\begin{lemma}[Uncovered indices are barriers]
\label{lem:uncovered-barrier}
Any uncovered index $i\in[n]$ is a barrier.
\end{lemma}
By Lemma~\ref{lem:uncovered-barrier}, to establish that an index $i\in[n]$ is a barrier, it suffices to show that it is uncovered. Consequently, verifying the barrier-free property reduces to the following probabilistic question: how likely is it that a long interval has at least one uncovered index?
The main difficulty is that the event ``index $i$ is uncovered'' depends on {all} neighborhoods $\cB_j$. As a result, uncovered-index events are generally dependent, even though the coordinates $\alpha_j$ themselves are independent. We address this issue by separating coverage due to nearby indices from coverage due to remote indices.
\begin{definition}[Local and remote coverage] 
Given $D>0$, an index $i\in\mathbb{Z}$ is said to be \emph{$D$-locally uncovered} if $i\not\in \cB_j$ for every $j\in \mathbb Z$ such that $|i-j|\leq D$. An index $j\in\mathbb{Z}$ is called a \emph{$D$-remote coverer} of $i\in\mathbb{Z}$ if $|i-j|>D$ and $i\in \cB_j$.
\end{definition} 

\begin{figure}[t]
\centering

\begin{tikzpicture}[
    x=0.82cm,
    y=0.82cm,
    >=Latex,
    nodev/.style={
        circle,
        draw,
        minimum size=5.8mm,
        inner sep=0pt,
        font=\scriptsize
    },
    plain/.style={
        nodev,
        fill=white
    },
    covered/.style={
        nodev,
        fill=gray!18
    },
    source/.style={
        nodev,
        fill=gray!60,
        text=white
    },
    special/.style={
        nodev,
        fill=white,
        line width=0.9pt
    },
    covbox/.style={
        rounded corners=3pt,
        draw=black!65,
        dashed,
        line width=0.4pt,
        fill=gray!08
    },
    remotebox/.style={
        rounded corners=3pt,
        draw=black!70,
        dashed,
        line width=0.4pt,
        fill=gray!14
    },
    intbox/.style={
        draw=black!60,
        fill=gray!06
    }
]

\fill[intbox] (4.55,-0.88) rectangle (10.45,-0.08);
\draw[dashed, line width=0.4pt] (4.55,-1.00) -- (4.55,0.55);
\draw[dashed, line width=0.4pt] (10.45,-1.00) -- (10.45,0.55);

\node[font=\small] at (7.5,-1.22) {interval $\cI$};

\draw[covbox] (2.65,-0.42) rectangle (5.35,0.42);   %
\draw[covbox] (4.65,-0.42) rectangle (7.35,0.42);   %
\draw[covbox] (6.65,-0.42) rectangle (9.35,0.42);   %
\draw[covbox] (8.65,-0.42) rectangle (11.35,0.42);  %
\draw[covbox] (10.65,-0.42) rectangle (13.35,0.42); %

\node[font=\scriptsize] at (4,0.72) {$R_4=2$};
\node[font=\scriptsize] at (6,0.72) {$R_6=2$};
\node[font=\scriptsize] at (8,0.72) {$R_8=2$};
\node[font=\scriptsize] at (10,0.72) {$R_{10}=2$};
\node[font=\scriptsize] at (12,0.72) {$R_{12}=2$};

\foreach \i in {1,...,14}{
    \pgfmathtruncatemacro{\j}{\i+1}
    \draw[black!70, line width=0.5pt] (\i,0) -- (\j,0);
}

\foreach \i in {1,2,3,13,14,15}{
    \node[plain] at (\i,0) {\(\scriptstyle \i\)};
}
\foreach \i in {5,7,9,11}{
    \node[covered] at (\i,0) {\(\scriptstyle \i\)};
}
\node[source] at (4,0) {\(\scriptstyle 4\)};
\node[source] at (6,0) {\(\scriptstyle 6\)};
\node[source] at (8,0) {\(\scriptstyle 8\)};
\node[source] at (10,0) {\(\scriptstyle 10\)};
\node[source] at (12,0) {\(\scriptstyle 12\)};

\node[font=\small] at (8,-1.85)
{(a) Every node in $I$ is covered locally.};

\begin{scope}[yshift=-2.8cm]

\fill[intbox] (4.55,-0.88) rectangle (10.45,-0.08);
\draw[dashed, line width=0.4pt] (4.55,-1.00) -- (4.55,0.55);
\draw[dashed, line width=0.4pt] (10.45,-1.00) -- (10.45,0.55);

\node[font=\small] at (7.0,-1.22) {interval $\cI$};

\draw[covbox] (2.65,-0.42) rectangle (5.35,0.42);   %
\draw[covbox] (4.65,-0.42) rectangle (7.35,0.42);   %
\draw[covbox] (6.65,-0.42) rectangle (9.35,0.42);   %
\draw[covbox] (10.65,-0.42) rectangle (13.35,0.42); %

\draw[remotebox] (9.65,-0.42) rectangle (15.35,0.42); %

\node[font=\scriptsize] at (4,0.72) {$R_4=2$};
\node[font=\scriptsize] at (6,0.72) {$R_6=2$};
\node[font=\scriptsize] at (8,0.72) {$R_8=2$};
\node[font=\scriptsize] at (12,0.72) {$R_{12}=2$};
\node[font=\scriptsize] at (14,0.72) {$R_{14}=5$};

\foreach \i in {1,...,14}{
    \pgfmathtruncatemacro{\j}{\i+1}
    \draw[black!70, line width=0.5pt] (\i,0) -- (\j,0);
}

\foreach \i in {1,2,3,13,15}{
    \node[plain] at (\i,0) {\(\scriptstyle \i\)};
}
\foreach \i in {5,7,9,11}{
    \node[covered] at (\i,0) {\(\scriptstyle \i\)};
}
\node[source] at (4,0)  {\(\scriptstyle 4\)};
\node[source] at (6,0)  {\(\scriptstyle 6\)};
\node[source] at (8,0)  {\(\scriptstyle 8\)};
\node[source] at (12,0) {\(\scriptstyle 12\)};
\node[source] at (14,0) {\(\scriptstyle 14\)};

\node[special] at (10,0) {\(\scriptstyle 10\)};

\draw[->, line width=0.5pt]
    (10,-1.48) -- (10,-0.36);
\node[font=\small] at (10,-1.70)
    {locally uncovered};

\draw[->, line width=0.5pt]
    (14,-1.10) -- (14,-0.36);
\node[font=\small] at (14,-1.32)
    {remote coverer};

\node[font=\small] at (8,-2.25)
{(b) A node in $\cI$ is not covered locally, but a remote coverer reaches it.};

\end{scope}
\end{tikzpicture}

\caption{Two ways in which an interval can contain no uncovered nodes, illustrated for $D=3$. In panel (a), every node in the interval is covered by a node within distance $D$. In panel (b), node $10$ is $D$-locally uncovered, but it is covered by the remote coverer $14$.}
\label{fig:local-remote-coverage}
\end{figure}
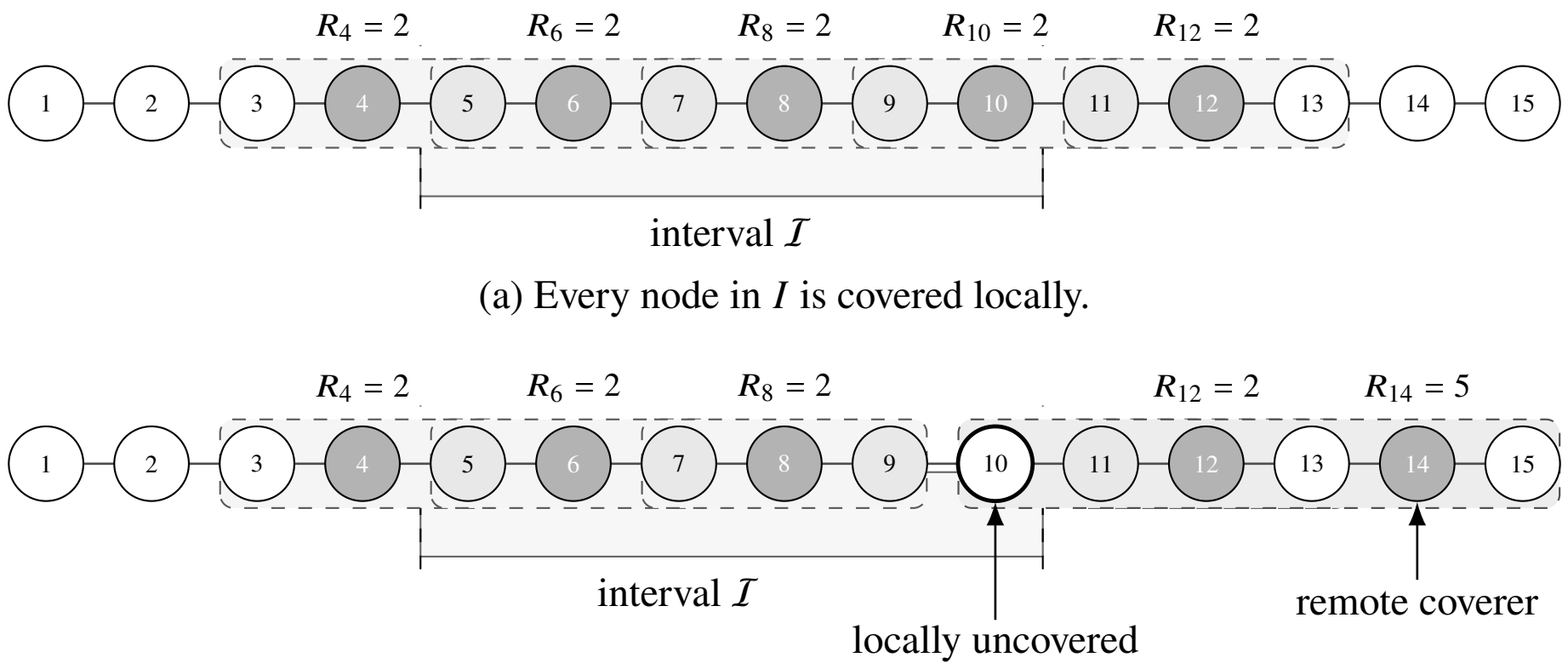

Figure~\ref{fig:local-remote-coverage} illustrates these two coverage mechanisms. The local--remote decomposition separates the dependence issue into two parts. First, $D$-local uncoveredness depends only on a block of coordinates around each candidate index. Thus, candidates separated by more than $2D$ have disjoint coordinate blocks, making their $D$-local-uncovered events independent under the product-measure assumption. Second, once sufficiently many locally uncovered candidates are identified, it remains to bound the probability that they are covered by remote indices. The next lemma formalizes the first observation.

\begin{lemma}[Probability and independence of local-uncovered events]
\label{lem:local-uncovered-independent}
Under the assumptions of Proposition~\ref{prop:explicit-product-barrier-free},
fix $D\in\mathbb{N}$ and, for each $i\in[n]$, let $A_i$
denote the event that $i$ is $D$-locally uncovered. 
Then, there exists a constant $p_D>0$, depending only on $\tau_b, g,\rho, q,p_0, M_q$ and independent of $n$, such that, for every $i\in[n]$,
\[
    \mathbb{P}_{\bm\alpha\sim\mu}(A_i)
    \geq p_D.
\]
Moreover, if $i,j\in[n]$ satisfy $|i-j|\geq 2D+1$, then $A_i$ and $A_j$ are independent.
\end{lemma}

\begin{lemma}[Many locally uncovered candidates]
\label{lem:many-local-uncovered}
Under the assumptions of
Lemma~\ref{lem:local-uncovered-independent}, fix $D\in\mathbb{N}$,
and let $c_1,\ldots,c_m\in[n]$ be deterministic candidate indices,
where $m\geq1$, whose pairwise distances are at least $2D+1$.
Let $X
    :=
    \left|
        \{k\in[m]:c_k\text{ is $D$-locally uncovered}\}
    \right|$
denote the number of locally uncovered candidates. Then
\[
    \mathbb{P}_{\bm\alpha\sim\mu}
    \left(X>\frac{p_Dm}{2}\right)
    \geq 1-
    \exp\left(-\frac{p_Dm}{8}\right).
\]
\end{lemma}

Choose candidates $c_1,\ldots,c_m\in\cI\subseteq[n]$ spaced
$2D+1$ apart. By Lemma~\ref{lem:many-local-uncovered}, a positive
fraction are $D$-locally uncovered with high probability for
sufficiently large $m$. However, these candidates may still have remote coverers
(see Figure~\ref{fig:local-remote-coverage}). Define
\[
    N_{\geq D}(j)
    :=\bigl|\{k\in[m]:D<|j-c_k|,\; c_k\in \cB_j\}\bigr|,
    \qquad j\in\mathbb Z
\]
as the number of candidates for which index $j$ is a $D$-remote coverer. The next lemma bounds $\sum_{j\in\mathbb Z}N_{\geq D}(j)$,
the total number of remote coverings counted with multiplicity,
by a small fraction of $m$ with high probability.

\begin{lemma}[Few remotely covered candidates]
\label{lem:few-remote-coverings}
Under the assumptions of Proposition~\ref{prop:explicit-product-barrier-free}, there exist an integer \(D\geq1\) and a constant $0<p_R\leq p_D$, depending only on $\tau_b,g,\rho,q,p_0,M_q$ and independent of $n$, such that the following holds: for every collection of candidate indices \(c_1,\ldots,c_m\) separated by at least \(2D+1\),

$$
    \mathbb{P}_{\bm\alpha\sim\mu}
    \left(
        \sum_{j\in\mathbb{Z}} N_{\geq D}(j)
        \leq
        \frac{p_R m}{2}
    \right)
    \geq
    1-
    \exp\left(-\frac{3p_R m}{8}\right).
$$
\end{lemma}
By Lemmas~\ref{lem:many-local-uncovered} and~\ref{lem:few-remote-coverings}, with high probability, locally uncovered candidates outnumber remote coverings. Thus, a barrier exists by Lemma~\ref{lem:uncovered-barrier}, proving Proposition~\ref{prop:explicit-product-barrier-free}.

\subsection{Compact SOCP Representation}\label{subsec::SOCP}

We now present an SOCP-representable formulation of
$\conv(\cP_{M,b})$ that agrees with the restricted DP
\eqref{eq:DP-recurrence-restricted}, including the case $\zero\in\Gamma$. Let $M\ge b$ be an integer;
as before, the reset condition is vacuous when $M>n$.
Define the acyclic graph $G=(\cV,\cA)$ by
\begin{align*}
    \cV&=\{0\}\cup[1:n-b+1]\cup\{n+1\},\\
    \cA&=\{(0,j):j\in \cV\setminus\{0\},\ j\le M\}\\
     &\quad\cup
       \{(i,j):i\in[1:n-b+1],\ j\in \cV,\ b\le j-i\le b+M-1\}.
\end{align*}
In the graph, vertex $0$ is a source vertex and $n+1$ is a terminal vertex.
Each internal vertex $i$ marks the start of a chosen resetting
sequence. Arc $(0,j)$ represents the interval $[1:j-1]$ before the first chosen resetting sequence. 
Arc $(i,j)$ with $i\ge1$ represents the block
$[i:j-1]$, consisting of $b$ zeros followed by a reset-free
interval of length at most $M-1$. 
For integers $0\le\ell\le k$ and a positive definite matrix
$\bm\Sigma\in\R^{k\times k}$, define
\begin{align*}
\mathcal X_{\ell,k}(\bm\Sigma)
    =\Big\{(t,\bm x,\bm z)&\in\R_+\times\R^k\times\{0,1\}^k:
       t\ge\bm x^\top\bm\Sigma\bm x,\quad
       \bm x\circ(\one-\bm z)=\zero,\\
    &z_i=0\quad\forall i\in[\ell],\qquad
       \sum_{i=j}^{j+b-1}z_i\ge1
       \quad\forall j\in[\ell+1:k-b+1]\Big\}.
\end{align*}
Constraints indexed by an empty interval are vacuous. For $k=0$,
we interpret $\mathcal X_{0,0}$ as the set consisting of $t\ge0$
and empty vectors $\bm x,\bm z$.
For each arc $(i,j)$, define $I_{ij}:=[\max\{1,i\}:j-1]$ and let $\bm Q_{I_{ij},I_{ij}}$ denote the submatrix of $\bm Q$ indexed by $I_{ij}$.
Let $\mathcal X^{ij}$ consist of the points
$(t,\bm x,\bm z)\in\cP$ with $x_\ell=z_\ell=0$ for $\ell\notin I_{ij}$,
$\bm z\in\cZ_{M,b}$, and
\begin{align*}
    (t,\bm x_{I_{0j}},\bm z_{I_{0j}})
       \in\mathcal X_{0,j-1}(\bm Q_{I_{0j},I_{0j}})\ \text{if }i=0,\qquad 
    (t,\bm x_{I_{ij}},\bm z_{I_{ij}})
       \in\mathcal X_{b,j-i}(\bm Q_{I_{ij},I_{ij}})\ \text{if }i\ge1.
\end{align*}
We discard every arc $(i,j)$ for which
$\mathcal X^{ij}=\emptyset$. We now show how, inspired by disjunctive programming representations of unions of convex sets \cite{ceria1999convex}, a description of $\conv(\mathcal X^{ij})$ can be used to construct an extended representation of $\conv(\cP_{M,b})$.  Specifically, let $f^{ij}$ be the convex indicator function of $\conv\left(\mathcal X^{ij}\right)$, that is, 
$$f^{ij}(t,\bm{x},\bm{z})=\begin{cases}0 &\text{ if }(t,\bm{x},\bm{z})\in \conv\left(\mathcal X^{ij}\right)\\
\infty &\text{otherwise.} \end{cases}$$
We defer to Proposition~\ref{prop:arc-support-hull} the discussion of a mechanism to compute functions $f^{ij}$.
Introduce flow variables $r^{ij}$ and local variables
$(t^{ij},\bm{x^{ij}},\bm{z^{ij}})\in\R\times\R^n\times\R^n$
for $(i,j)\in \cA$. Additional variables
$r^\star,t^\star\in\R$ and $\bm x^\star\in\R^n$ represent the
full-support $\bm z=\bm 1$. Consider

\begin{subequations}\label{eq:arc-convexification}
\begin{align}
    &\bm x=\sum_{(i,j)\in \cA}\bm{x^{ij}}+\bm x^\star,
       \qquad
       \bm z=\sum_{(i,j)\in \cA}\bm{z^{ij}}+\one r^\star,
       \label{eq:arc-linking}\\
    &t=\sum_{(i,j)\in \cA}t^{ij}+t^\star, \label{eq:arc-linking2}\\
    &r^{ij}f^{ij}\left(\frac{t^{ij}}{r^{ij}},\frac{\bm{x^{ij}}}{r^{ij}},\frac{\bm{z^{ij}}}{r^{ij}}\right)\leq 0
       &&\forall(i,j)\in \cA,\label{eq:arc-local-hulls}\displaybreak[4]\\
        &(\bm x^\star)^\top\bm Q\bm x^\star\le t^\star r^\star,
       \qquad t^\star\ge0,\quad r^\star\ge0,
       \label{eq:arc-full-support}\\
    &\sum_{(0,j)\in \cA}r^{0j}=1-r^\star,
       \qquad
       \sum_{(i,n+1)\in \cA}r^{i,n+1}=1-r^\star,
       \label{eq:arc-source-terminal}\\
    &\sum_{(i,k)\in \cA}r^{ik}-\sum_{(k,j)\in \cA}r^{kj}=0
       &&\forall k\in \cV\setminus\{0,n+1\},\\  
    &r^\star=0 &&\text{if }\zero\notin\Gamma,\\
    &r^\star\geq 0, \qquad r^{ij}\ge0 &&\forall(i,j)\in \cA,\label{eq:arc-ending}
\end{align}
\end{subequations}

where the function in \eqref{eq:arc-local-hulls} is the closure of the perspective function of $f^{ij}$.

\begin{proposition}\label{prop:arc-convexification}
The projection of~\eqref{eq:arc-convexification} onto
$(t,\bm x,\bm z)$ is $\conv(\cP_{M,b})$. 
\end{proposition}

\begin{proposition}\label{prop:arc-support-hull}
	Fix an arc $(i,j)$ with $\mathcal X^{ij}\ne\emptyset$, and let
	$\bm v^1,\ldots,\bm v^N$ list all its admissible support
	vectors, that is, all $\bm v\in\{0,1\}^n$ such that
	$(0,\zero,\bm v)\in\mathcal X^{ij}$.
	Then $f^{ij}(t,\bm x,\bm z)\leq 0$ if and only if there exist
	$\lambda_h,t^h\ge0$ and $\bm x^h\in\R^n$, $h\in[N]$, satisfying
	\begin{equation}\label{eq:arc-support-hull}
		\begin{aligned}
			&\sum_{h=1}^N\lambda_h=1,\qquad t\ge\sum_{h=1}^N t^h,\qquad \bm x=\sum_{h=1}^N\bm x^h,\qquad \bm z=\sum_{h=1}^N\lambda_h\bm v^h,\\
			&\bm x^h\circ(\one-\bm v^h)=\zero,
			h\in[N],\qquad (\bm x^h)^\top\bm Q\bm x^h\le t^h\lambda_h,
			h\in[N].
		\end{aligned}
	\end{equation}
\end{proposition}
\begin{proof}
	The disjunctive formulation of \citet{ceria1999convex}, applied to
	the quadratic epigraph for each fixed support $\bm v^h$, gives
	\eqref{eq:arc-support-hull}.\qed
\end{proof}

Let $\cY$ be obtained from~\eqref{eq:arc-convexification} by replacing each constraint~\eqref{eq:arc-local-hulls} with~\eqref{eq:arc-support-hull}, normalized by $\sum_h\lambda_h=r^{ij}$, and eliminating fixed-zero coordinates. The following theorem summarizes our main result: $\cY$ is an exact extended SOCP representation of $\conv(\cP_{M,b})$ and, under Theorem~\ref{thm:abstract-main}, yields a polynomial-size probabilistically exact convexification of $\cP$.

\begin{theorem}[Extended SOCP representation]\label{thm:compact-socp}
The projection of $\cY$ onto $(t,\bm x,\bm z)$ equals
$\conv(\cP_{M,b})$. Moreover, $\cY$ admits an SOCP description
in which the number of scalar variables, the number of linear
constraints, and the total second-order-cone dimension are
each bounded by $O\!\left(nM\Phi_b(M)\right)$
\begin{equation}\label{eq:direct-SOCP-size}
 O\!\left(nM\Phi_b(M)\right)
 =\begin{cases}
 O(nM^2),&b=1,\\
 O(nM\vartheta_b^M),&b\ge2.
 \end{cases}
\end{equation}
Consequently, under the assumptions of Theorem~\ref{thm:abstract-main} and with $M$ chosen as in~\eqref{eq:L-M}, $\conv(\cP_{M,b})$ is a
$(\mu,\Gamma,\delta)$-probabilistically exact convexification
of $\cP$ with an exact extended SOCP representation $\cY$ satisfying the size bounds
in~\eqref{eq:ideal-DD-size}.
\end{theorem}

\section*{Declaration of AI Use}
The authors independently developed the concept of probabilistically exact convexification and the initial proof strategy for banded structures in 2025, without the use of AI. They subsequently used ChatGPT 6 Pro to assist with scrutinizing mathematical arguments, checking proofs, and language polishing. All AI-assisted material was carefully reviewed and revised by the authors, who assume full responsibility for the mathematical correctness, content, and conclusions of this paper.

\section*{Acknowledgements}
Salar Fattahi is supported, in part, by grants CCF-2337776 and DMS-2152776 from the National Science Foundation and grant N000142612074 from the Office of Naval Research. 
Andr\'es G\'omez is supported, in part, by grant N00014-26-1-2117 from the Office of Naval Research and grant 2346058 from the National Science Foundation.

\bibliographystyle{plainnat}
\bibliography{ref}

\appendix
\section{Proofs}
\subsection{Proof of Lemma~\ref{lem:restricted-support-characterization}}
    We first show that both \(\conv(\cX)\) and
\(\conv(\cP )\) are closed. Let
\(\cC\) denote either \(\cX\) or \(\cP \). To prove that $\conv(\cC)$ is closed, we show that it can be written as the epigraph of a lower-semicontinuous function. By \cite[Theorem 1.6]{rockafellar1998variational}, this implies that $\cC$ is closed. 
Define $\varphi_{\cC}(\bm x,\bm z)
    :=
    \inf\{t:(t,\bm x,\bm z)\in\cC\},$
with the convention that \(\varphi_{\cC}(\bm x,\bm z)=+\infty\) when $\{t:(t,\bm x,\bm z)\in\cC\}$ is empty. Since \(\cC\) is upward-closed, $\cC
    =
    \left\{
        (t,\bm x,\bm z):
        t\geq\varphi_{\cC}(\bm x,\bm z)
    \right\},$ which is the epigraph of $\varphi(\bm x,\bm z)$. Since this epigraph is closed, 
we conclude that \(\varphi_{\cC}\) is lower semicontinuous \cite[Theorem 1.6]{rockafellar1998variational}.
Moreover, since \(\cC\subseteq\cP \),
\begin{align}
    \varphi_{\cC}(\bm x,\bm z)
    \geq \bm x^\top \bm Q\bm x
    \geq \lambda_{\min}\|\bm x\|_2^2.
\end{align}
Therefore, $\varphi$ is closed and $\frac{\varphi_{\cC}(\bm x,\bm z)}{\|(\bm x,\bm z)\|}\to+\infty$ whenever $\|(\bm x,\bm z)\|\to+\infty$ (recall that
\(\bm z\in\{0,1\}^n\) is bounded). Therefore, the conditions of \cite[Corollary 3.47]{rockafellar1998variational} are satisfied, which implies that the convex envelope of $\varphi_{\cC}$,
denoted by \(\conv\varphi_{\cC}\), is proper, lower semicontinuous, and coercive,
and its defining infimum over convex combinations is attained (see \cite[Proposition 2.31]{rockafellar1998variational}).
Therefore,
\begin{align}
    \conv(\cC)
    =
    \left\{
        (t,\bm x,\bm z):
        t\geq(\conv\varphi_{\cC})(\bm x,\bm z)
    \right\}.
\end{align}
Thus, \(\conv(\cC)\) is the epigraph of a lower-semicontinuous
function and is therefore closed.

Now, set $\widehat{\cX}:=\conv(\cX)$ and $\widehat{\cP}:=\conv(\cP )$.
As established above, both sets are closed and convex. Moreover, they are
upward-closed and
\(\widehat{\cX}\subseteq\widehat{\cP}\).
Suppose, toward a contradiction, that this inclusion is strict. Choose $(\bar t,\bm{\bar x},\bm{\bar z})
    \in \widehat{\cP}\setminus\widehat{\cX}.$
Because \(\widehat{\cX}\) is nonempty, closed, and convex, the strong
separating hyperplane theorem \cite[Corollary~11.4.2]{rockafellar1970convex} yields
\((\gamma,\bm\alpha,\bm\beta)\neq(0,\bm 0,\bm 0)\) and \(\eta\in\mathbb{R}\) such that
\begin{align}
    \gamma\bar t+\bm\alpha^\top\bm{\bar x}+\bm\beta^\top\bm{\bar z}
    <
    \eta
    \leq
    \gamma t+\bm\alpha^\top \bm x+\bm\beta^\top\bm z,
    \qquad
    \forall (t,\bm x,\bm z)\in\widehat{\cX}.
    \label{eq:strict-separation}
\end{align}
We must have \(\gamma\geq0\). Indeed, if \(\gamma<0\), then for any
\((t,\bm x,\bm z)\in\widehat{\cX}\), upward-closedness implies that
\((t+s,\bm x,\bm z)\in\widehat{\cX}\) for every \(s\geq0\). Consequently, $\gamma(t+s)+\bm\alpha^\top\bm x+\bm\beta^\top\bm z\to-\infty$ as $s\to+\infty$,
contradicting~\eqref{eq:strict-separation}.
If \(\gamma=0\), define $\Delta
    :=
    \eta-\bm\alpha^\top\bm{\bar x}-\bm\beta^\top\bm{\bar z}
    >0.$
Since \(\bar t\geq0\), choose \(\gamma_+>0\) sufficiently small that
\(\gamma_+\bar t<\Delta\). Then
\begin{align}
    \gamma_+\bar t+\bm\alpha^\top\bm{\bar x}+\bm\beta^\top\bm{\bar z}
    <
    \eta
    \leq
    \gamma_+ t+\bm\alpha^\top\bm x+\bm\beta^\top\bm z,
    \qquad
    \forall(t,\bm x,\bm z)\in\widehat{\cX},
\end{align}
where the second inequality follows from $\gamma_+>\gamma$ and \(t\geq0\). Thus, replacing
\(\gamma\) by \(\gamma_+\), we may assume without loss of generality
that \(\gamma>0\).
Dividing the strict separation inequality by \(\gamma\) gives
\begin{align}
    \bar t
    +\left(\frac{\bm\alpha}{\gamma}\right)^\top\bm{\bar x}
    +\left(\frac{\bm\beta}{\gamma}\right)^\top\bm{\bar z}
    <
    h_{\widehat{\cX}}
    \left(
        \frac{\bm\alpha}{\gamma},
        \frac{\bm\beta}{\gamma}
    \right).
\end{align}
On the other hand, since
\((\bar t,\bm{\bar x},\bm{\bar z})\in\widehat{\cP}\),
\begin{align}
    h_{\widehat{\cP}}
    \left(
        \frac{\bm\alpha}{\gamma},
        \frac{\bm\beta}{\gamma}
    \right)
    \leq
    \bar t
    +\left(\frac{\bm\alpha}{\gamma}\right)^\top\bm{\bar x}
    +\left(\frac{\bm\beta}{\gamma}\right)^\top\bm{\bar z}.
\end{align}
Therefore,
\begin{align}
    h_{\widehat{\cP}}
    \left(
        \frac{\bm\alpha}{\gamma},
        \frac{\bm\beta}{\gamma}
    \right)
    <
    h_{\widehat{\cX}}
    \left(
        \frac{\bm\alpha}{\gamma},
        \frac{\bm\beta}{\gamma}
    \right).
\end{align}
Recalling that $\widehat{\cX}:=\conv(\cX)$ and $\widehat{\cP}:=\conv(\cP )$, we have $h_{\widehat{\cX}}=h_{\cX}$ and $h_{\widehat{\cP}}=h_{\cP }$, which contradicts the assumed equality of the restricted lower-support
functions. Hence,
\[
    \conv(\cX)=\conv(\cP ).
\]
The converse follows immediately from $h_{\widehat{\cX}}=h_{\cX}$ and $h_{\widehat{\cP}}=h_{\cP }$.\qed

\subsection{Proof of Lemma~\ref{lem:barrier-forces-zeros}}
The proof relies on the following auxiliary lemma, which is a special case of \cite[Lemma 5]{bhathena2026solving}.
\begin{lemma}[inverse decay; Lemma 5 of \cite{bhathena2026solving}]\label{lem:inverse-decay}
Let $\cT\subseteq [n]$ be any subset of the rows/columns of $\bm Q$, and let $\pi: \cT\to [n]$ be the canonical indexing map that assigns to each $i\in \cT$ its corresponding row/column position within the submatrix $\bm Q_{\cT,\cT}$. For all indices $i,j\in \cT$,
\begin{equation}\label{eq:graph-inverse-decay}
\left|\left(\bm Q_{\cT,\cT}^{-1}\right)_{\pi(i),\pi(j)}\right|
\le \frac{2}{\lambdamin}\rho^{{|i-j|}},
\end{equation}
\end{lemma}

\begin{proof}[Proof of Lemma~\ref{lem:barrier-forces-zeros}]
    Fix \(\beta\in[\beta_0,\infty)^n\) and an optimal solution
\((t^\star, \bm x^\star,\bm z^\star)\) of~\eqref{eq:MIQP}, and let
\(\cT:=\{j\in[n]:z_j^\star=1\}\). The vector \(\bm x^\star\) minimizes \(\bm x^\top\bm Q\bm x+\bm\alpha^\top\bm x\) subject
to \(\bm x_{T^c}=0\). Hence
\[
    \bm x_\cT^\star=-\frac{1}{2}\bm Q_{\cT,\cT}^{-1}\bm\alpha_\cT,
    \qquad
    \bm x_{\cT^c}^\star=0.
\]
Invoking Lemma~\ref{lem:inverse-decay}, for every \(j\in T\),
\begin{align*}
    |x_j^\star|
    &\le \frac{1}{2}
    \sum_{q\in\cT}
    \left|
    \left(\bm Q_{\cT,\cT}^{-1}\right)_{\pi(j),\pi(q)}
    \right||\alpha_q|\le \frac{1}{\lambdamin}
    \sum_{q\in\cT}\rho^{|j-q|}|\alpha_q|
    \le \frac{1}{\lambdamin}W_j(\bm\alpha).
\end{align*}
Due to the definition of the $W_i$ (Definition~\ref{def:local-weight}) for every $k\ge0$ with $i+k\le n$,
\begin{align}\label{eq:W-ik}
    W_{i+k}(\bm\alpha)\leq \rho^{-k}W_i(\bm\alpha).
\end{align}
Set $\xi:=\sqrt{\frac{\beta_0}{2\lambdamax}}$.
For every \(k\in\{0,\ldots,b-1\}\), using
\eqref{eq:W-ik} and the definition of \(\tau_b\),
\[
    |x_{i+k}^\star|
    \le
    \frac{1}{\lambdamin}W_{i+k}(\bm\alpha)
    \le
    \frac{\rho^{-k}}{\lambdamin}W_i(\bm\alpha)
    \le
    \rho^{b-1-k}\xi
    \le
    \xi.
\]
Let \(\cA:=\cT\cap[i:i+b-1]\). We show that \(\cA=\varnothing\). Suppose, toward a contradiction, that \(\cA\neq\varnothing\). Define \(\bm{\widetilde x}\) by setting the coordinates of \(\bm x^\star\) indexed by \(\cA\) to zero and leaving all other coordinates unchanged. Then \(\bm{\widetilde x}\) is feasible for the support \(\cT\setminus \cA\).
For any support \(\cS\subseteq[n]\), write
\[
    f_{\bm\alpha,\bm\beta}(\bm x,\cS)
    :=
    \bm x^\top\bm Q\bm x+\bm\alpha^\top\bm x+\beta(\cS),
    \qquad
    \beta(\cS):=\sum_{j\in \cS}\beta_j.
\]
Since \(\bm x^\star\) minimizes the objective over the support
\(\cT\), its first-order optimality conditions give
\((2\bm Q\bm x^\star+\bm\alpha)_\cT=0\). By expanding $f_{\bm\alpha,\bm\beta}(\bm x^\star,\cT)$ and $f_{\bm\alpha,\bm\beta}(\bm{\widetilde x},\cT\setminus \cA)$ over $\cA$ and $\cT\setminus\cA$, we obtain
\begin{align}
    &f_{\bm\alpha,\bm\beta}(\bm x^\star,\cT)
    -f_{\bm\alpha,\bm\beta}(\bm{\widetilde x},\cT\setminus \cA)
    \nonumber\\
    &\qquad=
    (\bm x_\cA^\star)^\top\bm Q_{\cA\cA}\bm x_\cA^\star
    +2(\bm x_\cA^\star)^\top\bm Q_{\cA,\cT\setminus\cA}
    \bm x_{\cT\setminus\cA}^\star
    +\bm\alpha_\cA^\top\bm x_\cA^\star
    +\beta(\cA)
    \nonumber\\
    &\qquad=
    2(\bm x_\cA^\star)^\top(\bm Q\bm x^\star)_\cA
    -(\bm x_\cA^\star)^\top\bm Q_{\cA\cA}\bm x_\cA^\star
    +\bm\alpha_\cA^\top\bm x_\cA^\star
    +\beta(\cA)
    \nonumber\\
    &\qquad=
    \beta(\cA)
    -(\bm x_\cA^\star)^\top\bm Q_{\cA\cA}\bm x_\cA^\star.
    \label{eq:removal-identity}
\end{align}
On the other hand, \(\bm Q_{\cA\cA}\preceq\lambdamax\bm I\) and
\(|x_j^\star|\le\xi\) for every \(j\in \cA\), so
\[
    (\bm x_\cA^\star)^\top\bm Q_{\cA\cA}\bm x_\cA^\star
    \le
    \lambdamax\|\bm x_\cA^\star\|_2^2
    \le
    \lambdamax |\cA|\xi^2
    =
    \frac{|\cA|\beta_0}{2}
    \le
    \frac{\beta(\cA)}{2}.
\]
Substituting this bound into~\eqref{eq:removal-identity} gives
\[
    f_{\bm\alpha,\bm\beta}(\bm x^\star,\cT)
    -f_{\bm\alpha,\bm\beta}(\bm{\widetilde x},\cT\setminus \cA)
    \ge
    \frac{\beta(\cA)}{2}
    >0,
\]
where the strict inequality follows from \(\cA\neq\varnothing\) and
\(\beta_j\ge\beta_0>0\). Thus, removing the indices in \(\cA\) strictly
improves the objective, contradicting the optimality of
\((\bm x^\star,\bm z^\star)\). Hence \(\cA=\varnothing\), and therefore
\(z_i^\star=\cdots=z_{i+b-1}^\star=0\).\qed
\end{proof}

\subsection{Proof of Theorem~\ref{thm:abstract-main}}

Define
\[
    L
    :=
    M-b+1
    =
    \left\lceil
        \frac{1}{c_{\rm v}}
        \log\left(\frac{C_{\rm v}n}{\delta}\right)
    \right\rceil.
\]
If \(M>n\), then there is no block of \(M\) consecutive coordinates.
Under the convention that the defining condition of \(\cZ_{M,b}\) is
then vacuous, we have
\(\widehat{\cZ}_{M,b}=\{0,1\}^n\) and hence
\(\cP_{M,b}=\cP\). The conclusion is therefore immediate. We may
thus assume throughout the remainder of the proof that \(M\le n\).
For each \(s\in[n-M+1]\), let \(\cB_s:=[s:s+M-1]\) be the corresponding
block of \(M\) consecutive coordinates. The indices at which a
resetting sequence of length \(b\) can begin while remaining entirely
within \(B_s\) form the interval
\[
    \cI_s
    :=
    [s:s+M-b]
    =
    [s:s+L-1].
\]
In particular, \(|\cI_s|=L\) and
\(\cI_s\subseteq[n-b+1]\). By the
\((C_{\rm v},c_{\rm v})\) barrier-free property,
\[
    \Pp_{\bm\alpha\sim\mu}
    \left(
        W_i(\bm\alpha)>\tau_b
        \text{ for every }i\in\cI_s
    \right)
    \le
    C_{\rm v}e^{-c_{\rm v}L}
    \le
    \frac{\delta}{n},
\]
where the second inequality follows from the definition of \(L\).
Define the event
\[
    G
    :=
    \left\{
        \bm\alpha\in\R^n\;\mid\;
        \text{for every }s\in[n-M+1],
        \text{ there exists }i\in\cI_s
        \text{ such that }W_i(\bm\alpha)\le\tau_b
    \right\}.
\]
Since there are \(n-M+1\le n\) blocks of length \(M\), a union bound
gives
\[
    \Pp_{\bm\alpha\sim\mu}(G^c)
    \le
    (n-M+1)\frac{\delta}{n}
    \le
    \delta.
\]
where $G^c$ is the complement of $G$.
Thus, \(\Pp_{\bm\alpha\sim\mu}(G)\ge 1-\delta\).
Now, fix \(\bm\alpha\in G\). We show that
\[
    h_{\cP_{M,b}}(\bm\alpha,\bm\beta)
    =
    h_{\cP}(\bm\alpha,\bm\beta)
    \qquad
    \text{for every }\bm\beta\in\Gamma.
\]
First, consider any \(\bm\beta\in[\beta_0,\infty)^n\), and let
\((t^\star,\bm x^\star,\bm z^\star)\) be an optimal solution over \(\cP\).
For every \(s\in[n-M+1]\), the definition of \(G\) yields
an index \(i\in\cI_s\) such that \(W_i(\bm\alpha)\le\tau_b\).
Lemma~\ref{lem:barrier-forces-zeros} then implies $z_i^\star=z_{i+1}^\star=\cdots=z_{i+b-1}^\star=0.$
Since \(i\in\cI_s\), this resetting sequence is contained entirely
within \(\cB_s\). As this holds for every \(s\in[n-M+1]\), every block
of \(M\) consecutive coordinates of \(\bm z^\star\) contains a resetting
sequence of length \(b\). Hence
\(\bm z^\star\in\cZ_{M,b}\subseteq\widehat{\cZ}_{M,b}\), and therefore
\((t^\star,\bm x^\star,\bm z^\star)\in\cP_{M,b}\). This implies \(h_{\cP_{M,b}}(\bm\alpha,\bm\beta)\le h_{\cP}(\bm\alpha,\bm\beta)\).
Because \(\cP_{M,b}\subseteq\cP\), we have
\(h_{\cP_{M,b}}(\bm\alpha,\bm\beta)\ge h_{\cP}(\bm\alpha,\bm\beta)\). Thus,
\[
    h_{\cP_{M,b}}(\bm\alpha,\bm\beta)
    =
    h_{\cP}(\bm\alpha,\bm\beta)
    \qquad
    \forall\bm\beta\in\Gamma\cap[\beta_0,\infty)^n.
\]
It remains to consider the possible case \(\bm \beta=\zero\). Suppose
that \(\zero\in\Gamma\). The all-one support imposes no restrictions
on \(\bm x\), and it belongs to
\(\widehat{\cZ}_{M,b}\) by definition. In particular, setting
\[
    \bm{\bar x}:=-\frac{1}{2}\bm Q^{-1}\bm\alpha,
    \qquad
    \bm{\bar z}:=\one,
    \qquad
    \bar t:=\bm{\bar x}^\top\bm Q\bm{\bar x}
\]
gives a point \((\bar t,\bm{\bar x},\bm{\bar z})\in\cP_{M,b}\). Moreover,
\(\bm{\bar x}\) is the unique minimizer of
\(\bm x^\top\bm Q\bm x+\bm\alpha^\top\bm x\) over \(\R^n\). Since the all-one support
also belongs to \(\cP\), this point is optimal over \(\cP\) when
\(\bm\beta=\zero\). Therefore,
\[
    h_{\cP_{M,b}}(\bm\alpha,\zero)
    =
    h_{\cP}(\bm\alpha,\zero).
\]
Under the standing assumption
\(\Gamma\subseteq\{\zero\}\cup[\beta_0,\infty)^n\), for every
\(\bm\alpha\in G\), we have 
\(h_{\cX}(\bm\alpha,\bm\beta)=h_{\cP }(\bm\alpha,\bm\beta)\) for all $\bm\beta\in \Gamma$. Hence
\[
    \mathbb{P}_{\bm\alpha\sim\mu}
        \left(
            h_{\cX}(\bm\alpha,\bm\beta)
            =
            h_{\cP }(\bm\alpha,\bm\beta)
            \quad\text{for all }\bm\beta\in\Gamma
        \right)
    \ge
    1-\delta.
\]
Finally, \(\cP_{M,b}\) is nonempty and upward-closed in \(t\). Therefore, \(\cP_{M,b}\) is a
\((\mu,\Gamma,\delta)\)-probabilistically exact approximation of
\(\cP\) according to
Definition~\ref{def:inner-approx-convex}. By the same definition,
\(\conv(\cP_{M,b})\) is a \((\mu,\Gamma,\delta)\)-probabilistically exact convexification.\qed

\subsection{Proof of Lemma~\ref{lem:uncovered-barrier}}
By definition of $R_j$ and $R_j\leq |i-j|$, we have
$|\alpha_j|\le a_gg^{|i-j|}$.  Therefore
\[
W_i(\bm\alpha) = \sum_{j=1}^n \rho^{|i-j|}|\alpha_j| \leq \sum_{j\in\Z}\rho^{|i-j|}|\alpha_j|
\le a_g\sum_{d\in\Z}(\rho g)^{|d|}
=a_g\frac{1+\rho g}{1-\rho g}=\tau_b.
\]\qed

\subsection{Proof of Lemma~\ref{lem:local-uncovered-independent}}

We first derive a uniform tail bound for the coverage radii. By the
definition of $R_j$, for every $j\in[n]$ and $d\in\mathbb{Z}_+$,
\[
    \{R_j>d\}
    =
    \{|\alpha_j|>a_gg^d\}.
\]
Define $C_q:=M_qa_g^{-q}$ and $\psi_q:=g^{-q}\in(0,1)$.
Markov's inequality and
\eqref{eq:uniform-product-conditions} give
\begin{equation}
\label{eq:radius-lower-bound-local}
    \mathbb{P}(R_j>d)
    \leq
    \frac{\mathbb{E}|\alpha_j|^q}
         {\left(a_g g^{d}\right)^{q}}
    \leq
    M_qa_g^{-q} g^{-qd}
    =
    C_q\psi_q^d. \quad \implies\quad \mathbb{P}(R_j\leq d)
    \geq
    1-C_q\psi_q^d.
\end{equation}
On the other hand, since $g^d\geq1$,
\[
    \{|\alpha_j|\leq a_g\}
    \subseteq
    \{|\alpha_j|\leq a_gg^d\}
    =
    \{R_j\leq d\}.
\]
This implies that $\mathbb{P}(R_j\leq d)\geq p_0$ since, by the first condition of \eqref{eq:uniform-product-conditions}, we have $\mathbb{P}(|\alpha_j|\leq a_g)\geq p_0$.
Combining the two bounds yields
\begin{equation}
\label{eq:radius-prob-lower}
    \mathbb{P}(R_j\leq d)
    \geq
    \max\{p_0,1-C_q\psi_q^d\}.
\end{equation}
Now fix $i\in[n]$. The event $A_i$ can be written as
\begin{align*}
        A_i
    &=\{i\not\in \cB_i\}
    \cap
    \bigcap_{d=1}^{D}
    \left(\{i\not\in \cB_{i-d}\}
    \cap
    \{i\not\in \cB_{i+d}\}\right)\\
    &=
    \{R_i=0\}
    \cap
    \bigcap_{d=1}^{D}
    (\{R_{i-d}\leq d\}
    \cap
    \{R_{i+d}\leq d\}),
\end{align*}
Each $R_j$ is a function only of $\alpha_j$. Since $\mu$ is a product
measure, the radii associated with distinct indices in $[n]$ are
independent. Moreover,
\[
    \mathbb{P}(R_i=0)
    =
    \mathbb{P}(|\alpha_i|\leq a_g)
    \geq p_0.
\]
Using independence together with
\eqref{eq:radius-prob-lower}, and noting that for each
$d=1,\ldots,D$ there are at most two relevant indices at distance
$d$ from $i$, we obtain
\[
    \mathbb{P}(A_i)
    \geq p_D:=
    p_0
    \prod_{d=1}^{D}
    \max\left\{
        p_0,\,
        1-C_q\psi_q^d
    \right\}^2.
\]
Finally, $A_i$ depends only on the coordinates indexed by $[i-D:i+D]\cap[n],$
while $A_j$ depends only on those indexed by $[j-D:j+D]\cap[n].$
If $|i-j|\geq 2D+1$, these two index sets are disjoint. Since $\mu$
is a product measure, the corresponding coordinate blocks are
independent. Therefore, $A_i$ and $A_j$ are independent.\qed

\subsection{Proof of Lemma~\ref{lem:many-local-uncovered}}

For each $k\in[m]$, recall that $A_{c_k}$ denotes the event that $c_k$
is $D$-locally uncovered. The separation of the candidates together with Lemma~\ref{lem:local-uncovered-independent} implies that the events
$A_{c_1},\ldots,A_{c_m}$ are mutually independent.
By Lemma~\ref{lem:local-uncovered-independent},
$\mathbb{P}(A_{c_k})\geq p_D$ for every $k\in[m]$. Thus
$X=\sum_{k=1}^{m}\mathbf{1}_{A_{c_k}}$ is a sum of independent
Bernoulli variables with $\mathbb{E}X\geq p_Dm$.
The Chernoff bound gives
\[
    \mathbb{P}\left(X>\frac{p_Dm}{2}\right)
    \geq
    1-\mathbb{P}\left(X\leq\frac{\mathbb{E}X}{2}\right)
    \geq
    1-\exp\left(-\frac{\mathbb{E}X}{8}\right)
    \geq
    1-\exp\left(-\frac{p_Dm}{8}\right).
\]
\qed

\subsection{Proof of Lemma~\ref{lem:few-remote-coverings}}

For brevity, write $S:=\sum_{j\in\mathbb{Z}}N_{\geq D}(j)$.
Recall that $\alpha_j=0$ for $j\notin[n]$, and hence $R_j=N_{\geq D}(j)=0$ for such
indices. Thus, the sum defining $S$ is effectively finite. Define 
$$
    p_R
    :=
    p_0
    \prod_{d=1}^{\infty}
    \max\left\{
        p_0,\,
        1-C_q\psi_q^d
    \right\}^2>0,
$$
where $C_q:=M_qa_g^{-q}$ and $\psi_q:=g^{-q}\in(0,1)$, as defined before.
We first show that $p_R>0$. Let $p_d:=\max\{p_0,1-C_q\psi_q^d\}$ for $d\geq1$. Since
$0<p_0\leq p_d\leq1$ and
\[
    \sum_{d=1}^{\infty}|p_d-1|
    \leq
    C_q\sum_{d=1}^{\infty}\psi_q^d
    =
    \frac{C_q\psi_q}{1-\psi_q}
    <\infty,
\]
the standard infinite-product convergence criterion
\cite[Theorem~1]{trench1999conditional} implies that
$\prod_{d=1}^{\infty}p_d$ converges to a strictly positive limit.
Consequently, $p_R=p_0\prod_{d=1}^{\infty}p_d^2>0$.

Set $s_D:=2D+1$. Since the candidate indices $c_1,\ldots,c_m$ are
separated by at least $s_D$, a radius $R_j$ can cover at most
$1+2R_j/s_D$ candidate indices. Consequently,
$N_{\geq D}(j)\leq (1+2R_j/s_D)\mathbf{1}_{\{R_j>D\}}$.
For fixed candidate indices, $N_{\geq D}(j)$ is a deterministic function of $R_j$
alone. Since $R_j$ depends only on $\alpha_j$ and $\mu$ is a product
measure, the variables $(N_{\geq D}(j))_{j\in\mathbb{Z}}$ are independent.
Therefore,
\[
    \log\mathbb{E}e^S
    =
    \sum_{j\in\mathbb{Z}}\log\mathbb{E}e^{N_{\geq D}(j)}\leq \sum_{j\in\mathbb{Z}}\left(\mathbb{E}e^{N_{\geq D}(j)}-1\right)\leq \sum_{j\in\mathbb{Z}}\mathbb{E}\!\left[N_{\geq D}(j)e^{N_{\geq D}(j)}\right],
\]
where the equality follows from the independence of $(N_{\geq D}(j))_{j\in\mathbb{Z}}$; the first inequality follows from $\log u\leq u-1$ for every $u>0$; and the second inequality follows from $e^u-1\leq ue^u$ for $u\geq0$. 
By the definition of $N_{\geq D}(j)$, we have
$N_{\geq D}(j)=\sum_{k=1}^m
\mathbf{1}_{\{D<|j-c_k|<R_j\}}$. Moreover,
$N_{\geq D}(j)\leq1+2R_j/s_D$. Therefore, for each $j\in \mathbb{Z}$, we have
\[
N_{\geq D}(j)e^{N_{\geq D}(j)} = \sum_{k=1}^m \mathbf{1}_{\{D<|j-c_k|<R_j\}} e^{N_{\geq D}(j)}\leq e\sum_{k=1}^m \mathbf{1}_{\{D<|j-c_k|<R_j\}} e^{2R_j/s_D}.
\]
Hence,
\[
\begin{aligned}
    \log\mathbb{E}e^S
    &\leq e\sum_{j\in \mathbb{Z}}\sum_{k=1}^m \mathbb{E}\!\left[\mathbf{1}_{\{D<|j-c_k|<R_j\}} e^{2R_j/s_D}\right]\\
    &\leq
    e\sum_{k=1}^m
    \sum_{\substack{j\in\mathbb{Z}\\|j-c_k|>D}}
    \mathbb{E}\!\left[
        \mathbf{1}_{\{R_j>|j-c_k|\}}e^{2R_j/s_D}
    \right]                                                     \\
    &\leq
    2em
    \sum_{d=D+1}^{\infty}
    \sup_{j\in[n]}
    \mathbb{E}\!\left[
        \mathbf{1}_{\{R_j>d\}}e^{2R_j/s_D}
    \right],
\end{aligned}
\]
where the factor $2$ accounts for the two possible indices at each
distance $d$ from a candidate index.
Define
$u_D:=\psi_q e^{2/s_D}$. Recall that $\psi_q\in (0,1)$. Therefore, for $D$ sufficiently large, $u_D<1$. For every
$d\geq D+1$,
\[
\begin{aligned}
    \mathbb{E}\!\left[
        \mathbf{1}_{\{R_j>d\}} e^{2R_j/s_D}
    \right]
    &=
    \sum_{r=d+1}^{\infty}
    e^{2r/s_D}\mathbb{P}(R_j=r)                                      \\
    &\leq
    \sum_{r=d+1}^{\infty}
    e^{2r/s_D}\mathbb{P}(R_j> r-1)                                  \\
    &\leq
    \frac{C_q}{\psi_q}
    \sum_{r=d+1}^{\infty}u_D^r\\
    &
    =
    \frac{C_q}{\psi_q}
    \frac{u_D^{d+1}}{1-u_D},
\end{aligned}
\]
where the last inequality follows from the fact that, for every
$j\in[n]$ and $r\in\mathbb{Z}_+$, $\mathbb{P}(R_j>r-1)\leq C_q\psi_q^{r-1}.$
Substituting this estimate above yields
\[
    \log\mathbb{E}e^S
    \leq
    m\varepsilon_D,
    \quad \text{where}\quad
    \varepsilon_D
    :=
    \frac{2eC_q}{\psi_q}
    \frac{u_D^{D+2}}{(1-u_D)^2}.
\]
Since $u_D=\psi_q e^{2/(2D+1)}\to\psi_q<1$ as $D\to\infty$, we have
$\varepsilon_D\to0$ as $D\to\infty$. We may therefore choose a fixed $D$, depending only on $g,\rho,q,p_0,M_q$, such that
$\varepsilon_D\leq p_R/8$.
Finally, Markov's inequality gives
\[
\begin{aligned}
    \mathbb{P}
    \left(
        S>\frac{p_R m}{2}
    \right)
    \leq
    \exp\left(-\frac{p_R m}{2}\right)
    \mathbb{E}e^S             \leq
    \exp\left[
        -\left(\frac{p_R}{2}-\varepsilon_D\right)m
    \right]
    \leq
    \exp\left(-\frac{3p_R m}{8}\right).
\end{aligned}
\]
Equivalently,
\[
    \mathbb{P}
    \left(
        \sum_{j\in\mathbb{Z}}N_{\geq D}(j)
        \leq
        \frac{p_R m}{2}
    \right)
    \geq
    1-\exp\left(-\frac{3p_R m}{8}\right),
\]
which proves the claim.\qed

\subsection{Proof of Proposition~\ref{prop:explicit-product-barrier-free}}
Let $0<p_R\leq p_D$ and $D\geq1$ be as in Lemma~\ref{lem:few-remote-coverings}. Fix a nonempty interval $\cI\subseteq[n-b+1]$ and select $m=\lceil|\cI|/(2D+1)\rceil$ candidates in $\cI$ separated by at least \(2D+1\). By Lemmas~\ref{lem:local-uncovered-independent} and~\ref{lem:many-local-uncovered} and noting $p_R\leq p_D$, more than $p_R m/2$ candidates are $D$-locally uncovered, except on an event of probability at most $e^{-p_R m/8}$. Meanwhile, Lemma~\ref{lem:few-remote-coverings} ensures that the total number of remote coverings of these candidates, counted with multiplicity, is at most $p_R m/2$, except on an event of probability at most $e^{-3p_R m/8}$.
When both conclusions hold, at least one locally uncovered candidate has no remote coverer. This candidate is globally uncovered and hence is a barrier by Lemma~\ref{lem:uncovered-barrier}. A union bound therefore gives
\[
    \mathbb{P}_{\bm\alpha\sim\mu}
    \left(
        W_i(\bm\alpha)>\tau_b
        \text{ for every }i\in\cI
    \right)
    \leq
    e^{-p_R m/8}+e^{-3p_R m/8}
    \leq
    2\exp\left(-\frac{p_R|\cI|}{8(2D+1)}\right).
\]
Thus, $\mu$ satisfies the barrier-free property with $C_{\rm v}=2$ and $c_{\rm v}=p_R/[8(2D+1)]>0$. These constants are independent of $n$ and $\cI$.\qed

\subsection{Proof of Proposition~\ref{prop:arc-convexification}}
	Consider the minimization of an arbitrary linear function of the form
	\begin{align}
		\min\;&\bm{\alpha^\top x}+\bm{\beta^\top z}+t\text{ s.t. }\eqref{eq:arc-convexification}.\label{eq:minimization Hull1}
	\end{align}
	We show that minimization \eqref{eq:minimization Hull1} is equivalent to the restricted DP \eqref{eq:DP-recurrence-restricted}, and thus is equivalent to minimizing the same function over $\cP_{M,b}$. First, using constraints \eqref{eq:arc-linking}-\eqref{eq:arc-linking2}, we can rewrite \eqref{eq:minimization Hull1} as
		{\small\begin{align}
		\min\;&\bm{\alpha^\top}\left(\sum_{(i,j)\in \cA}\bm{x^{ij}}+\bm x^\star\right)+\bm{\beta^\top}\left(\sum_{(i,j)\in \cA}\bm{z^{ij}}+\one r^\star\right)+\sum_{(i,j)\in \cA}t^{ij}+t^\star\text{ s.t. }\eqref{eq:arc-local-hulls}-\eqref{eq:arc-ending}.\label{eq:minimization Hull2}
	\end{align}}
Now observe that if variables $\bm{r}$ are fixed, problem \eqref{eq:minimization Hull2} is separable across variable $(t^{ij},\bm{x^{ij}},\bm{z^{ij}})$ and $(t^{\star},\bm{x^{\star}},\bm{1})$. If $r^{ij}=0$, then constraint \eqref{eq:arc-local-hulls} forces $\bm{z^{ij}}=\bm{x^{ij}}=\bm{0}$ and $t^{ij}=0$ in optimal solutions; otherwise, if $i\geq 1$, then the cost associated with variables $(t^{ij},\bm{x^{ij}},\bm{z^{ij}})$ is
\begin{align}
	&\min_{(t^{ij},\bm{x^{ij}},\bm{z^{ij}})}\bm{\alpha^\top x^{ij}}+\bm{\beta^\top z^{ij}}+t^{ij}\text{  s.t.  }r^{ij}f^{ij}\left(\frac{t^{ij}}{r^{ij}},\frac{\bm{x^{ij}}}{r^{ij}},\frac{\bm{z^{ij}}}{r^{ij}}\right)\leq 0\notag \\
	\Leftrightarrow\;& \min_{(\bar t^{ij},\bm{\bar x^{ij}},\bm{\bar z^{ij}})}r^{ij}\left(\bm{\alpha^\top \bar x^{ij}}+\bm{\beta^\top \bar z^{ij}}+\bar t^{ij}\right)\text{  s.t.  }f^{ij}(\bar t^{ij},\bm{\bar x^{ij}},\bm{\bar z^{ij}})\leq 0,\label{eq:eq:minimization Hull3}
\end{align}
where the equivalence is obtained by doing the change of variables $(\bar t^{ij},\bm{\bar x^{ij}},\bm{\bar z^{ij}})=\frac{1}{r^{ij}}( t^{ij},\bm{ x^{ij}},\bm{ z^{ij}})$. From the definition of $f^{ij}$ as a convex indicator function, we find that
\begin{align}
	\eqref{eq:eq:minimization Hull3}\Leftrightarrow\;&\min_{(\bar t^{ij},\bm{\bar x^{ij}},\bm{\bar z^{ij}})}r^{ij}\left(\bm{\alpha^\top \bar x^{ij}}+\bm{\beta^\top \bar z^{ij}}+\bar t^{ij}\right)\text{  s.t.  }(\bar t^{ij},\bm{\bar x^{ij}},\bm{\bar z^{ij}})\in \conv(\mathcal{X}^{ij})\notag\\
	=&\min_{(\bar t^{ij},\bm{\bar x^{ij}},\bm{\bar z^{ij}})}r^{ij}\left(\bm{\alpha^\top \bar x^{ij}}+\bm{\beta^\top \bar z^{ij}}+\bar t^{ij}\right)\text{  s.t.  }(\bar t^{ij},\bm{\bar x^{ij}},\bm{\bar z^{ij}})\in \mathcal{X}^{ij}\notag\\
	=&r^{ij}\widehat q_{[i+b:j-1]},\label{eq:simplifiedCost1}
\end{align}
where $\widehat q_{[i+b:j-1]}$ is the restricted DP cost of the
reset-free interval following the $b$ zeros at $[i:i+b-1]$, as defined in Section~\ref{subsec::DP}.
 Observe that cost \eqref{eq:simplifiedCost1} also holds when
$r^{ij}=0$. Through identical arguments, the cost associated with
a source arc $(0,j)$ is $r^{0j}\widehat q_{[1:j-1]}$, since its
associated interval $[1:j-1]$ must be reset-free. The cost associated with
$(t^{\star},\bm{x^{\star}},\bm{1})$ is
$r^\star \left(-\frac{1}{4 }\bm{\alpha^\top Q^{-1}\alpha}
+\bm{\beta^\top 1}\right)$, which is the cost of minimizing the
quadratic program assuming a full support multiplied by $r^\star$.
Thus, projecting out all auxiliary variables except $\bm{r}$ from
\eqref{eq:minimization Hull2}, we find the equivalent form
\begin{align*}
    \min\;&\sum_{\substack{(i,j)\in \cA\\i\ge1}}
             r^{ij}\widehat q_{[i+b:j-1]}
       +\sum_{(0,j)\in \cA}r^{0j}\widehat q_{[1:j-1]}\\
       &\quad+r^\star\left(-\frac14\bm{\alpha^\top Q^{-1}\alpha}
                              +\bm{\beta^\top 1}\right)\\
    \text{s.t. }&\eqref{eq:arc-source-terminal}-\eqref{eq:arc-ending}.
\end{align*}
This optimization problem corresponds precisely to a shortest path problem over the state transition graph of the restricted DP \eqref{eq:DP-recurrence-restricted}, and thus their optimal objective values coincide.

We next verify that both $\conv(\cP_{M,b})$ and the projection
of~\eqref{eq:arc-convexification} onto $(t,\bm x,\bm z)$ are closed.
The set $\cP_{M,b}\subseteq \cP$ is upward-closed. Hence, by Lemma~\ref{lem:restricted-support-characterization}, $\conv(\cP_{M,b})$ is closed.
For the projection, consider any convergent sequence of projected
points and choose a feasible lift of each point. The flow variables
lie in $[0,1]$, the nonnegative local epigraph variables are bounded
by $t$, and the perspective constraints together with $\bm Q\succ0$
bound all remaining local variables. The lifts therefore have a
convergent subsequence, whose limit is feasible because the extended
feasible set is closed. Thus, the limiting projected point also
belongs to the projection, proving its closedness.

Both sets are therefore convex, upward-closed, and contained
in $\{(t,\bm x,\bm z):t\ge0\}$. Since their optimal values agree for
every objective $t+\bm\alpha^\top\bm x+\bm\beta^\top\bm z$, the
separation argument in the proof of
Lemma~\ref{lem:restricted-support-characterization} yields the
claimed convex-hull equality.\qed

\subsection{Proof of Theorem~\ref{thm:compact-socp}}
\label{app:restricted-dp}
We first provide an upper and lower bounds on the number of reset-free patterns, which will be useful in establishing a bound on the size of $\cY$.

\begin{lemma}[Number of reset-free patterns]\label{lem:word-count} 
Let $N_b(k)$ be the number of binary support vectors of length $k$ containing no $b$ consecutive zeros, with $N_b(0)=1$. Set $\vartheta_1=1$ and, for $b\ge2$, let $\vartheta_b\in(1,2)$ be the unique positive root of $$x^b=x^{b-1}+\cdots+x+1.$$ There are constants $c_b,C_b>0$, depending only on $b$, such that 
\begin{equation}\label{eq:word-count} 
c_b\vartheta_b^k\le N_b(k)\le C_b\vartheta_b^k \quad \forall k\geq 0, \qquad \sum_{k=0}^{M-1}N_b(k)\le C_b\Phi_b(M), 
\end{equation} 
where 
\[
    \Phi_b(M)
    :=
    \begin{cases}
        M, & b=1,\\
        \vartheta_b^M, & b\ge2.
    \end{cases}
\]
\end{lemma}
\begin{proof}
If $b=1$, the only reset-free vector of each length is the all-ones vector, including the empty vector when $k=0$. Thus $N_1(k)=1$ and $\sum_{k=0}^{M-1}N_1(k)=M=\Phi_1(M)$, so both bounds hold with $c_1=C_1=1$. Henceforth assume $b\ge2$.

For $0\le k<b$, every binary vector of length $k$ is reset-free,
and hence $N_b(k)=2^k$. We first derive a recurrence for $N_b(k)$.
Fix $k\ge b$. For a binary vector $\bm z$ of length $k$, let $t(\bm z)$ denote
its number of trailing zeros, with $t(\bm z)=0$ when $z_k=1$.
Partition the reset-free vectors according to
$t(\bm z)$. Since such a vector cannot contain $b$ consecutive zeros,
$t(\bm z)\in\{0,\ldots,b-1\}$. If $t(\bm z)=j$, its last
$j+1$ entries are a one followed by $j$ zeros. Deleting these
entries leaves a reset-free vector of length $k-j-1$. Conversely,
appending a one followed by $j$ zeros to any reset-free vector of
length $k-j-1$ produces a reset-free vector with exactly $j$ trailing
zeros. This bijection gives
\begin{equation}\label{eq:reset-free-count-recursion}
    N_b(k)=\sum_{j=0}^{b-1}N_b(k-j-1)
          =\sum_{j=1}^{b}N_b(k-j),\qquad k\ge b.
\end{equation}
The polynomial $p_b(x):=x^b-x^{b-1}-\cdots-x-1$ has at most
one positive root, since $p_b(x)/x^b=1-\sum_{j=1}^b x^{-j}$
is strictly increasing for $x>0$. Moreover,
\[
    p_b(1)=1-b<0,\qquad
    p_b(2)=2^b-\sum_{j=0}^{b-1}2^j=1>0.
\]
Thus its unique positive root $\vartheta_b$ lies in $(1,2)$.
Define the positive constants
\[
    c_b:=\min_{0\le r<b}\frac{N_b(r)}{\vartheta_b^r},
    \qquad
    D_b:=\max_{0\le r<b}\frac{N_b(r)}{\vartheta_b^r}.
\]
We claim that $c_b\vartheta_b^k\le N_b(k)\le D_b\vartheta_b^k$
for every $k\ge0$. The claim holds for $0\le k<b$ by construction.
Suppose it holds for all indices smaller than some $k\ge b$.
Using~\eqref{eq:reset-free-count-recursion} and the defining equation
for $\vartheta_b$, we obtain
\begin{align*}
    N_b(k)
    &=\sum_{j=1}^{b}N_b(k-j)
      \le D_b\sum_{j=1}^{b}\vartheta_b^{k-j}\\
    &=D_b\vartheta_b^{k-b}
      (\vartheta_b^{b-1}+\cdots+\vartheta_b+1)
      =D_b\vartheta_b^k.
\end{align*}
The same induction hypothesis gives
$N_b(k)\ge c_b\sum_{j=1}^{b}\vartheta_b^{k-j}
=c_b\vartheta_b^k$, proving the claim.
Finally, since $\vartheta_b>1$, the cumulative count satisfies
\begin{align*}
    \sum_{k=0}^{M-1}N_b(k)
    &\le D_b\sum_{k=0}^{M-1}\vartheta_b^k
     =D_b\frac{\vartheta_b^M-1}{\vartheta_b-1}\\
    &\le\frac{D_b}{\vartheta_b-1}\vartheta_b^M.
\end{align*}
Taking $C_b=D_b\max\{1,(\vartheta_b-1)^{-1}\}$ makes the same
constant valid in both upper bounds in~\eqref{eq:word-count}.\qed
\end{proof}

\begin{proof}[Proof of Theorem~\ref{thm:compact-socp}]
For $r^{ij}>0$, dividing the local variables by $r^{ij}$
recovers~\eqref{eq:arc-support-hull} with normalization
$\sum_h\lambda_h=1$. When $r^{ij}=0$, nonnegativity forces
$\lambda_h=0$ for every $h$, and $\bm Q\succ0$ then forces
$\bm x^h=\zero$. The resulting local points are exactly
$\bm x^{ij}=\bm z^{ij}=\zero$, $t^{ij}\ge0$, as required
by the closed perspective. Thus, \eqref{eq:arc-support-hull} with the normalization
$\sum_h\lambda_h=r^{ij}$ is an equivalent SOCP representation of
\eqref{eq:arc-local-hulls}. The projection equality then follows
from Proposition~\ref{prop:arc-convexification}.

To bound the size, fix an arc $(i,j)\in A$ and let $k<M$ be the length of the reset-free interval in the definition of $\cX^{ij}$; explicitly, $k=j-1$ if $i=0$ and $k=j-i-b$ otherwise. Since all support coordinates outside this interval are fixed to zero, each admissible support of $\cX^{ij}$ is determined by a reset-free binary vector of length $k$. Thus, $\cX^{ij}$ has at most $N_b(k)$ admissible support patterns. After omitting coordinates fixed to zero, each pattern requires at most $k$ coordinates of $\bm x^h$, the two nonnegative scalars $\lambda_h,t^h$, and a second-order cone of dimension at most $k+2$. The local vectors $\bm x^{ij}$ and $\bm z^{ij}$ each retain at most $k$ coordinates. Together with $t^{ij}$, $r^{ij}$, and the aggregation and normalization constraints, these contribute an additional $O(k+1)$ scalar variables and linear constraints. Hence the local SOCP formulation associated with \eqref{eq:arc-local-hulls} has $O((k+1)N_b(k))$ scalar variables, linear constraints, and total second-order-cone dimension.

For each $k$, there are at most $n+1$ such arcs, since the
starting vertex and $k$ determine the endpoint.
The global variables, flow balance equations, and
full-support block contribute only $O(n)$ additional terms.
Consequently, the size of $\cY$ is bounded by
\[
    O\!\left(n\sum_{k=0}^{M-1}(k+1)N_b(k)\right)
    =O_b\!\left(nM\Phi_b(M)\right),
\]
where we used $k+1\le M$ and
$\sum_{k=0}^{M-1}N_b(k)=O_b(\Phi_b(M))$ from
Lemma~\ref{lem:word-count}. This gives
\eqref{eq:direct-SOCP-size}.
Finally, Theorem~\ref{thm:abstract-main} gives probabilistic
exactness for the choice~\eqref{eq:L-M}. Substituting
$M\le b+c_{\rm v}^{-1}\log(C_{\rm v}n/\delta)$ into
\eqref{eq:direct-SOCP-size} yields~\eqref{eq:ideal-DD-size},
with $\gamma_b=\log(\vartheta_b)/c_{\rm v}$.\qed
\end{proof}

\section{Specific Distributions Satisfying the Barrier-Free Property}\label{app:explicit-dist}
\begin{corollary}[Explicit families satisfying the barrier-free property]
\label{lem:explicit-product-families}
The conditions of
Proposition~\ref{prop:explicit-product-barrier-free} hold in
each of the following cases, with the displayed choices of
$q$, $p_0$, and $M_q$.

\begin{enumerate}
    \item \textbf{Gaussian distributions.}
    Suppose that $\alpha_i\sim\mathcal{N}(0,\sigma_i^2)$ with
    $0<\sigma_i\leq\bar \sigma<\infty$ for every $i\in[n]$.
    One may take $q=1$, $p_0=2\Phi(a_g/\bar \sigma)-1$, $M_1=\bar \sigma\sqrt{\frac{2}{\pi}}$,
    where $\Phi$ is the standard normal distribution function.

    \item \textbf{Laplace distributions.}
    Suppose that $\alpha_i$ has density
    $f_i(x)=(2\sigma_i)^{-1}\exp(-|x|/\sigma_i)$, with
    $0<\sigma_i\leq\bar \sigma<\infty$ for every $i\in[n]$.
    One may take $q=1,
        p_0=1-\exp(-a_g/\bar \sigma),
        M_1=\bar \sigma.$

    \item \textbf{Bounded distributions with positive mass near zero.}
    Suppose that there exist constants $B>0$ and
    $\pi_0\in(0,1]$ such that $|\alpha_i|\leq B$ almost surely
    and $\mathbb{P}(|\alpha_i|\leq a_g)\geq\pi_0$ for every
    $i\in[n]$. One may take $q=1,
        p_0=\pi_0,
        M_1=B.$

    \item \textbf{Cauchy distributions.}
    Suppose that $\alpha_i\sim\operatorname{Cauchy}(0,\sigma_i)$
    with $0<\sigma_i\leq\bar \sigma<\infty$ for every $i\in[n]$.
    One may take $q=\frac12,
        p_0=\frac{2}{\pi}\arctan(a_g/\bar \sigma),
        M_{1/2}=\sqrt{2\bar \sigma}.$

    \item \textbf{Student-$t$ distributions.}
    Suppose that $\alpha_i=\sigma_iT_i$, where each $T_i$ has a
    standard Student-$t_\nu$ distribution for a fixed $\nu>0$,
    and $0<\sigma_i\leq\bar \sigma<\infty$ for every $i\in[n]$.
    For any fixed $q\in(0,\nu)$, one may take
    \[
        p_0=2F_\nu(a_g/\bar \sigma)-1,
        \qquad
        M_q=
        \bar \sigma^q\nu^{q/2}
        \frac{
            \Gamma\bigl((q+1)/2\bigr)
            \Gamma\bigl((\nu-q)/2\bigr)
        }{
            \sqrt{\pi}\,\Gamma(\nu/2)
        },
    \]
    where $F_\nu$ is the standard Student-$t_\nu$
    distribution function and $\Gamma$ is the gamma function.
\end{enumerate}

\end{corollary}

\begin{proof}
For each distribution with chosen parameters, it suffices to verify the local-mass and moment bounds in
\eqref{eq:uniform-product-conditions}.

We first start with the bounded family. The local-mass condition for this family is assumed,
and the moment bound $\mathbb{E}|\alpha_i|\leq B$ follows from
$|\alpha_i|\leq B$ almost surely.

We next consider the local-mass condition for the Gaussian,
Laplace, Cauchy, and Student-$t$ families. In each case,
$\alpha_i=\sigma_iZ_i$, where the standardized variables $Z_i$
have a common continuous distribution function $F$ that is
symmetric about zero. Since $\sigma_i\leq\bar \sigma$,
\[
    \mathbb{P}(|\alpha_i|\leq a_g)
    =
    2F(a_g/\sigma_i)-1
    \geq
    2F(a_g/\bar \sigma)-1.
\]
For the Gaussian and Student-$t$ families, this gives the choices of $p_0$ directly. For the standard Laplace
distribution, $2F(u)-1=1-e^{-u}$ for $u\geq0$, while for the
standard Cauchy distribution,
$2F(u)-1=(2/\pi)\arctan u$. These identities give the
remaining choices of $p_0$. Each is strictly positive because
$a_g/\bar \sigma>0$ and the corresponding density is positive in
a neighborhood of zero.

We next verify the moment bounds. For a standard normal
variable $Z_G$ and a standard Laplace variable $Z_L$,
direct integration gives
\[
    \mathbb{E}|Z_G|
    =
    \frac{2}{\sqrt{2\pi}}
    \int_0^\infty x e^{-x^2/2}\,dx
    =
    \sqrt{\frac{2}{\pi}},
    \qquad
    \mathbb{E}|Z_L|
    =
    \int_0^\infty x e^{-x}\,dx
    =
    1.
\]
Since $\mathbb{E}|\alpha_i|=\sigma_i\mathbb{E}|Z_i|$ and
$\sigma_i\leq\bar \sigma$, the asserted Gaussian and Laplace moment
bounds follow.
For a standard Student-$t_\nu$ variable $T$ and
$q\in(0,\nu)$, its density gives
\begin{align*}
    \mathbb{E}|T|^q
    &=
    \frac{2\Gamma((\nu+1)/2)}
         {\sqrt{\nu\pi}\,\Gamma(\nu/2)}
    \int_0^\infty
        x^q\left(1+\frac{x^2}{\nu}\right)^{-(\nu+1)/2}
    \,dx \\
    &=
    \frac{\nu^{q/2}\Gamma((\nu+1)/2)}
         {\sqrt{\pi}\,\Gamma(\nu/2)}
    \int_0^\infty
        u^{(q-1)/2}(1+u)^{-(\nu+1)/2}
    \,du \\
    &=
    \nu^{q/2}
    \frac{
        \Gamma((q+1)/2)\Gamma((\nu-q)/2)
    }{
        \sqrt{\pi}\,\Gamma(\nu/2)
    }.
\end{align*}
Here the second equality uses $u=x^2/\nu$, and the third
uses Euler's beta integral.
Multiplying by $\sigma_i^q\leq\bar \sigma^q$ proves the stated
Student-$t$ moment bound.

Finally, the standard Cauchy distribution is the
Student-$t$ distribution with $\nu=1$. Taking $q=1/2$
in the preceding identity yields
$\mathbb{E}|T|^{1/2}
=\Gamma(3/4)\Gamma(1/4)/\pi=\sqrt{2}$,
where the last equality follows from the gamma reflection
formula. Thus
$\mathbb{E}|\alpha_i|^{1/2}
=\sqrt{2\sigma_i}\leq\sqrt{2\bar \sigma}$, as claimed.\qed

\end{proof}

\end{document}